\documentclass[10pt]{article}
\usepackage{graphicx}
\usepackage{amssymb}
\usepackage{epstopdf}
\usepackage{color}
\usepackage{xcolor}
\usepackage{pstricks}
\usepackage{mdframed}
\usepackage{epic,eepic,epsfig,amssymb,amsmath,amsthm,graphics}
\usepackage{authblk}
\usepackage{amsmath}
\usepackage{hyperref}
\usepackage{amsfonts}
\usepackage{tikz}
\usepackage{tikz-3dplot}

\def\phi{{\varphi}}

\DeclareSymbolFont{AMSb}{U}{msb}{m}{n}
\DeclareMathSymbol{\N}{\mathbin}{AMSb}{"4E}
\DeclareMathSymbol{\Z}{\mathbin}{AMSb}{"5A}
\DeclareMathSymbol{\R}{\mathbin}{AMSb}{"52}
\DeclareMathSymbol{\Q}{\mathbin}{AMSb}{"51}
\DeclareMathSymbol{\I}{\mathbin}{AMSb}{"49}
\DeclareMathSymbol{\C}{\mathbin}{AMSb}{"43}

\DeclareMathOperator*{\argmax}{argmax}
\def\be{\begin{equation}}
\def\ee{\end{equation}}
\def\ber{\begin{eqnarray}}
\def\eer{\end{eqnarray}}

\def\beq{\begin{equation}}
\def\eeq{\end{equation}}

\newtheorem{theorem}{Theorem}[section]
\newtheorem{proposition}[theorem]{Proposition}
\newtheorem{lemma}[theorem]{Lemma}
\newtheorem{corollary}[theorem]{Corollary}
\theoremstyle{definition}
\newtheorem{remark}[theorem]{Remark}
\newtheorem{example}[theorem]{Example}

\newtheorem{definition}[theorem]{Definition}

\title{Stratified Monge-Kantorovich optimal transport problems
\footnote{A.M. is pleased to acknowledge the support of the  National Sciences and Engineering Research Council of Canada.}}

\author{  Mohammad Ali Ahmadpoor \thanks{School of Mathematics and Statistics,  Carleton University, Ottawa, ON,  Canada, mohammadaliahmadpoo@cmail.carleton.ca} \quad  and \quad Abbas Moameni\thanks{School of Mathematics and Statistics,  Carleton University, Ottawa, ON,  Canada, momeni@math.carleton.ca}
}
\date{}

\begin{document}
\maketitle
\begin{abstract}
Topology and uniqueness of   optimal plans   in the Monge-Kantorovich  optimization problem is one of the immensely important target of researches in this area. 
In this paper, we investigate Monge-Kantorovich optimal transport problems for which  the absolute continuity of marginal measures is relaxed.  For $X,Y\subseteq\mathbb{R}^{n+1}$ let $(X,\mathcal{B}_X,\mu)$ and $(Y,\mathcal{B}_Y,\nu)$ be two Borel probability spaces, $c:X\times Y\to\mathbb{R}$ be a cost function, and consider  the optimization problem
\begin{align*}\tag{MKP}\label{MKPEQ}
    \inf\left\{\int_{X\times Y} c(x,y)\,d\lambda\ :\ \lambda \in\Pi(\mu,\nu)  \right\}.
\end{align*}
 Inspired  by a  shape recognition problem in computer vision,  
in a seminal paper of Gangbo-McCann \cite{GANGBOMCCANN2} for the quadratic cost $c(x,y)=|x-y|^2$,  it is shown that  the uniqueness of Kantorovich solutions holds even though Monge solutions fail to
exist where the measure $\mu$ is allowed to charge the $n$-dimensional subsets of $X$. Motivated by this, our goal in this work is twofold. 
We first consider an optimal transport problem with multi-layers target space for the cost  $c(x,y)=h(x-y)$ with $h$ being  strictly convex. Namely,  we assume  
\[
    X=\overline{X}\times\{\overline{x}\},\quad\text{and}\quad Y=\bigcup_{k=1}^K \left(\overline{Y}_{k}\times \{\overline{y}_k\}\right),
\]
where  $\overline{X}, \overline{Y}_{k} \subseteq \mathbb{R}^{n}$ for $k\in \{1,\ldots,K\},$ $\overline{x}\in \mathbb{R}$,   and $\{\overline{y}_1,...,  \overline{y}_K\}\subseteq  \mathbb{R}$  is a set of  distinct numbers. 
We also assume that  the restriction of $\mu$ on $\overline{X}$ is absolutely continuous with respect to the $n$-dimensional  Lebesgue measure, but $\mu$ is singular with respect to the $(n+1)$-dimensional  Lebesgue measure.  When $K=1,$ this translates to the standard Monge-Kantorovich optimal transport problem for which the solution is unique and concentrates on a single map. We shall show that for $K\geq 2,$ the solution is still unique but it concentrates on the graph of several  maps. 
Secondly, for general closed subsets $X\subseteq \mathbb{R}^{n+1}$ we consider the case where 
 the first marginal $\mu$ is of the  form
\begin{align*}
    \int_X f(x)\,d\mu(x)=\int_X f(x)\alpha(x)\,d\mathcal{L}^{n+1}(x)+\int_{X_0} f(x_0)\, d S(x_0),\quad \forall f\in C_b(X)
\end{align*}
where $X_0 \subseteq X$ is  an $n$-dimensional manifold. Here, 
  $d\mathcal{L}^{n+1}$  stands for the $(n+1)$-dimensional Lebesgue measures on $\mathbb{R}^{n+1}$,  and $d S$ is a measure on $X_0$ which is absolutely continuous with respect to $d\mathcal{L}^{n}$ with respect to each coordinate chart on $X_0$.  This can be seen as a two-layers problem as  the measure $\mu$  charges both  $n$-dimensional and $n+1$-dimensional subsets.
\end{abstract}
\maketitle
\section{Introduction}
The Monge-Kantorovich optimal transport problem seeks a joint measure, denoted as $\lambda$, on the product space $X\times Y$ for marginal  Borel probability measures $\mu$ and $\nu$ defined on  domains $X$ and $Y$ in $\mathbb{R}^{n+1}$. The objective is to minimize the total transport cost 
$$C(\lambda) :=\int_{X\times Y} c(x, y) d\lambda(x, y),$$
where the cost function $c$ is  representing  the cost associated with moving a unit mass from $x\in X$ to $y\in Y$. The optimization is carried out within the convex set $\Pi(\mu,\nu)$ of joint measures that have $\mu$ and $\nu$ as marginals. If the optimal joint measure $\lambda$ is concentrated on the graph of a mapping $T:X\to Y$ that transports $\mu$ to $\nu$, then $\lambda$ can be expressed as $\lambda=(\text{id}\times T)\#\mu$. The map $T$, referred to as the optimal transport map, minimizes the total transport cost 
$$C(\lambda) :=\int_{X\times Y} c(x, T(x)) d\mu(x),$$
among all mappings that transport $\mu$ to $\nu$, a problem known as the Monge optimization (see \cite{VILLANI2009} as a thorough reference).

An illustrative example from Gangbo and McCann \cite{GANGBOMCCANN} highlights scenarios where such a map is applicable. These include cases where $\mu$ is absolutely continuous with respect to the Lebesgue measure on $\mathbb{R}^{n+1}$, and  $c(x,y)=h(x-y)$ with $h$ being a strictly convex function. Alternatively, when $\mu$ and $\nu$ have non-overlapping supports, and $c(x,y)=l(|x-y|)$ with $l$ being a strictly concave function. The optimal transport map provides insights into the magnitude and direction of mass movement in the vicinity of $x\in X$. Moreover, in \cite{FRY} for measures supported on the
domain boundaries $\partial X$ and $\partial Y$, if the cost is chosen to depend also on the relative orientation of the outward unit normals to these boundaries, Fry observed that the corresponding Monge-Kantorovich problem serves as a prototype for a shape recognition problem in computer vision that uses boundary matching as a form of comparison to identify objects. In a pioneering  paper Gangbo-McCann \cite{GANGBOMCCANN2} studied
$$\inf\left\{\int_{\partial X\times \partial Y} |x-y|^2 d\lambda(x, y)\ :\ \lambda\in\Pi(\tilde{\mu},\tilde{\nu})\right\},$$
where $\tilde{\mu}$ and $\tilde{\nu}$  are Borel probability measures on $\partial X$ and $\partial Y$, the boundaries of the bounded domains $X$ and $Y$ in $\mathbb{R}^{n+1}$. In contrast to the theory applicable to measures that are absolutely continuous concerning the Lebesgue measure on $\mathbb{R}^{n+1}$, it is observed that for $\mu$-almost every  $x\in\partial X$, there is not a unique destination $y\in\partial Y$ during transportation. Nevertheless, it has been demonstrated that the images of $x$ are almost surely arranged in a collinear manner and parallel the normal to $\partial X$ at $x$. In cases where either domain is strictly convex, they have concluded that the solution to the optimization problem is unique. When both domains exhibit uniform convexity, a regularity result has been established, indicating that the images of $x\in\partial X$ are consistently collinear. Additionally, both images vary with $x$ in a continuous and continuously invertible way.

Motivated by this, we explore the following problems. The first problem pertains to the strictly convex cost function $c(x,y) = h(x-y)$ in multi-layer spaces. Specifically, we delve into the two-marginal Monge-Kantorovich optimization problem
\begin{align}\label{STR-CONV-NORMAL-INTRO}
    \inf\left\{\int_{X\times Y}h(x-y)\,d\lambda\ :\ \lambda\in\Pi(\mu,\nu)\right\},
\end{align}
where  the function $h:\mathbb{R}^{n+1}\to\mathbb{R}^+$ is strictly convex, and the domains $X$ and $Y$ are differentiable $n$-dimensional manifolds in $\mathbb{R}^{n+1}$. In particular, for any compact subsets $\overline{X}, \overline{Y}_{k} \subseteq \mathbb{R}^{n}$ (for $k=1,\ldots,K$ with $K\in\mathbb{N}$), any element $\overline{x}$, and any set of distinct elements $\{\overline{y}_1,\ldots \overline{y}_K\}$ in $ \mathbb{R}$, we define $X$ and $Y$ as
\begin{align*}
    X=\overline{X}\times\{\overline{x}\},\quad\text{and}\quad Y=\bigcup_{k=1}^K \left(\overline{Y}_{k}\times \{\overline{y}_k\}\right).
\end{align*}
Our objective is to demonstrate that, owing to the geometric characteristics of the domain, the solution to problem \eqref{STR-CONV-NORMAL-INTRO} is unique but not induced by a single map. Indeed, in Theorem \ref{TH31}, it is shown that  the unique solution to \eqref{STR-CONV-NORMAL-INTRO} is concentrated on the union of the graphs of $K$ functions. In the case where $K=1$, the problem reduces to an optimization problem on $\mathbb{R}^n$, having a unique solution concentrated on the graph of a single optimal transport map \cite{GANGBOMCCANN}.

For the quadratic cost function $c(x, y) = |x - y|^2$, in the case where $\mu$ and $\nu$ have support on orthogonal subspaces of $\mathbb{R}^{n+1}$, it is known  that the solution to the minimization problem in \eqref{STR-CONV-NORMAL-INTRO} is not unique. Conversely, for any feasible measure $\lambda \in \Pi(\mu, \nu)$, it is proven that they all share the same cost. This is demonstrated by constructing a convex function $\phi$ whose subdifferential includes $\text{Spt}(\nu) \times \text{Spt}(\nu)$, as described in \cite[Remark 15]{MCCANN}. Specifically, one sets $\phi = 0$ on $\text{conv}(\text{Spt}(\mu))$ and $\phi= +\infty$ elsewhere. However, a more general scenario exists: in situations where the domains satisfy the condition $\langle\mathrm{n}(x) , \mathrm{n}(y)\rangle \neq 0$ for every $(x, y) \in X \times Y$, uniqueness is guaranteed, as outlined in Theorem \ref{TH32}. Here $\mathrm{n}(x)$ and $\mathrm{n}(y)$ denote the unit  normal vectors at the points $x$ and $y$.


The second problem we tackle involves cost functions $c$ satisfying the sub-twist condition \cite{AHMADKIMMCCANN}. For such functions, we consider a two-marginal optimal transport problem with the cost function $c:X\times Y\to\mathbb{R}$, where $X$ and $Y$ are two compact subsets of $\mathbb{R}^{n+1}$ with smooth boundaries, equipped with Borel probability measures $\mu$ and $\nu$, respectively. The measure $\mu$ takes the form
\begin{align*}
    \int_X f(x)\,d\mu(x)=\int_X f(x)\alpha(x)\,d\mathcal{L}^{n+1}(x)+\int_{X_0} f(x_0)\beta(x_0)\, d\mathcal{L}^{n}(x_0),\quad \forall f\in C_b(X),
\end{align*}   
Here, $\mathcal{L}^{n}$ denotes the Lebesgue measure on $\mathbb{R}^{n}$, and $\alpha, \beta$ are non-negative measurable functions. Additionally, $X_0$ is any sub-manifold of $X$ with dimension $n$  where  either lacks a boundary, or has a smooth $(n-1)$-dimensional boundary $\partial X_0$.
We shall establish the existence of an unique optimal plan that is not necessarily induced by a single map. We would like to remark that this uniqueness property is already  known if either $\alpha$ or $\beta$ is identically zero
\cite{AHMADKIMMCCANN}.



This paper is organized as follows. Section \ref{SECTIONPRELIMINARIES} provides a description of our setting and delves into the essential tools upon which our paper is based, with additional technical aspects covered in Section \ref{SECTIONAPPENDIX}. The problem \eqref{STR-CONV-NORMAL-INTRO} will be thoroughly examined in Section \ref{SECTIONMULTI-LAYER-STRCONV}, where we also explore the quadratic cost function and derive uniqueness in the most optimal way. Subsequently, Section \ref{SECTIONMULTI-LAYER-QUAD23} investigates the uniqueness of the optimal plan for the quadratic cost function across multi-layer domains. The handling of cost functions with the sub-twist property is discussed in Section \ref{SECTIONSUBTWIST}. Finally, Section \ref{SECTIONEXAMPLES} is dedicated to applications of our results.

\section{Preliminaries}\label{SECTIONPRELIMINARIES}

In this section, we introduce the current problem and its preliminaries, which will be needed in the subsequent sections. We commence with a concise overview of the optimal mass transportation problem. Let $X,Y\subseteq\mathbb{R}^{n+1}$, where $(X,\mathcal{B}_X,\mu)$ and $(Y,\mathcal{B}_Y,\nu)$ represent two Polish spaces, and $c:X\times Y\to\mathbb{R}$ is a cost function. We define the (measurable) map $T:X\to Y$ as pushing forward $\mu$ to $\nu$, denoted by $\nu=T\#\mu$, if
\begin{align*}
    \nu(B)=\mu(T^{-1}(B)),\quad\forall B\in\mathcal{B}_Y.
\end{align*}
Let $\mathcal{T}(\mu,\nu)$ denote the set of such maps. The Monge problem, as described in \cite{MONGE}, involves finding an optimal transport map $T\in\mathcal{T}(\mu,\nu)$ that minimizes the cost of transporting mass from $X$ to $Y$. In other words, $T$ is a solution to the following problem
\begin{align*}\tag{MP}\label{TAGMP}
    \inf\left\{\int_{X\times Y} c(x,T(x))\,d\mu(x)\ :\ T \in\mathcal{T}(\mu,\nu)  \right\}.
\end{align*}
The existence of a minimizer for \eqref{TAGMP} may not be expected under mild conditions or in a general setting, as it can depend on various factors. In \cite{KANTOROVICH}, a relaxation to problem \eqref{TAGMP} was introduced, marking an immensely important step in this area of research. It can be formulated as follows
\begin{align*}\tag{MKP}\label{TAGMKP}
    \inf\left\{\int_{X\times Y} c(x,y)\,d\lambda\ :\ \lambda \in\Pi(\mu,\nu)  \right\},
\end{align*}
where $\Pi(\mu,\nu)$ is the set of all probability measures whose first and second marginals on $X$ and $Y$ are equal to $\mu$ and $\nu$, respectively. In terms of notation,
\begin{align*}
    \lambda(A\times Y)=\mu(A),\quad \lambda(X\times B)=\nu(B),\quad \forall A\in\mathcal{B}_X,\ \forall B\in\mathcal{B}_Y.
\end{align*}


Moreover, the dual problem to \eqref{TAGMKP} can be stated as follows
\begin{align*}\tag{DMKP}\label{TAGDMKP}
    \sup\left\{\int_X \phi(x)\,d\mu +\int_Y \psi(y)\,d\nu \ :\ (\phi,\psi)\in L_1(X,\mu)\times L_1(Y,\nu),\ \phi+\psi\leq c \right\}.
\end{align*}
The problems \eqref{TAGMP} and \eqref{TAGMKP} have received extensive study in the literature. Indeed, due to their applications in a wide variety of branches of science such as economics, physics, fluid mechanics, and many other areas of mathematics (see \cite{AHMADKIMMCCANN,BIANCHINICARAVENNA,BUTTAZZODEPASCALEGORIGIORGI,CAFFARELLI,CHAMPIONDEPASCALE,CHIAPPORIMCCANNNESHEIM,FATHIFIGALLI,GANGBOTHESIS,GANGBOMCCANN,LEVIN,MOMENI-CHAR}, and \cite{VILLANI2009} for an overview), these problems have garnered the attention of researchers.

In this section, we present some fundamental theorems and definitions regarding the existence of solutions to the problems \eqref{TAGMKP} and \eqref{TAGDMKP} and some applicable results which can be viewed as powerful tools for dealing with the Monge-Kantorovich problems. It is a well-known fact that in our setting, i.e., in the situation in which the domains are compact and the cost function is continuous existence of the minimizer to the problem \eqref{TAGMKP} and maximizer to the dual problem \eqref{TAGDMKP} is guaranteed (for instance, see \cite{VILLANI2009}). We summarize these results in the following lemma.

\begin{lemma}\label{LEMMAEXISTENCE}
    Let $X$ and $Y$ be two compact subsets of $\mathbb{R}^{n+1}$ equipped with the probability measures $\mu$ and $\nu$, respectively and $c:X\times Y\to\mathbb{R}$ be a non-negative continuous function. Then there exists at least one minimizer $\lambda$ of the problem \eqref{TAGMKP}. Moreover, there exists continuous potential functions $\phi :X\to\mathbb{R}$ and $\psi :Y\to\mathbb{R}$ which maximize \eqref{TAGDMKP} and satisfy the following relations
    \begin{align*}
        &c(x,y)\geq \phi(x)+\psi(y),\quad \forall (x,y)\in X\times Y,\\
        &\phi(x)=\inf\left\{c(x,y)-\psi(y)\ :\ y\in Y\right\},\\
        &\psi(y)=\inf\left\{c(x,y)-\phi(x)\ :\ x\in X\right\}.
    \end{align*}
 Additionally, for any minimizer $\lambda$ of \eqref{TAGMKP}, we have that
 \begin{align*}
     \text{Spt}(\lambda)\subseteq \mathcal{S}:=\left\{ (x,y)\in X\times Y\ : \ c(x,y)=\phi(x)+\psi(y) \right\}.
 \end{align*}
    
\end{lemma}
One perspective on the potential functions here involves $c$-concave functions and the concept of $c$-superdifferential. Indeed, we have the following description. To explore these concepts further, one can refer to \cite{AGUIDE}.
\begin{definition}
\begin{enumerate}
    \item $c$-transform of a function $\psi:Y\to\mathbb{R}\cup\{\pm\infty\}$ is the function $\psi^c:X\to\mathbb{R}\cup\{-\infty\}$ defined by
    \begin{align}
        \psi^c(x)=\inf\left\{c(x,y)-\psi(y)\ :\ y\in Y\right\}. 
    \end{align}
    \item 
    A function $\phi:X\to\mathbb{R}\cup\{-\infty\}$ is said to be $c$-concave if there exists a function $\psi:Y\to\mathbb{R}\cup\{-\infty\}$ such that $\phi=\psi^c$.

    \item The $c$-superdifferntial of a $c$-concave function $\phi:X\to\mathbb{R}\cup\{-\infty\}$ is denoted by $\partial^c\phi$, and is the following subset of $X\times Y$,
    \begin{align*}
        \left\{(x,y)\in X\times Y\ :\ c(x,y)=\phi(x)+\phi^c(y)\right\}.
    \end{align*}

    \item The $c$-superdifferential of $\phi$ at a point $x\in X$ is therefore defined to be the $x$-intersection of the set $\partial^c\phi$, that is,
    \begin{align*}
        \partial^c\phi(x)=\left\{y\in Y\ :\ (x,y)\in\partial\phi\right\}.
    \end{align*}

    \end{enumerate}
\end{definition}
\begin{remark}
\begin{enumerate}
    \item 
The $c$-superdifferential of a $c$-concave function $\phi$, can be interpreted as follows. The pair $(x,y)\in X\times Y$ belongs to $\partial^c\phi$ if and only if
\begin{align*}
    c(x,y)-\phi(x)\leq c(x',y)-\phi(x'),\quad\forall x'\in X.
\end{align*}
Consequently, as a result, the $c$-superdifferential of a $c$-concave function is a $c$-cyclically monotone set, i.e., for any finite collection $\{(x_k,y_k)\}_{k=1}^n\subseteq\partial^c\phi$ and for any permutation mapping $\sigma :\{1,\ldots,n\}\to \{1,\ldots,n\}$, we have
        \begin{align*}
            \sum_{k=1}^n c(x_k,y_k)\leq \sum_{k=1}^n c(x_{k},y_{\sigma(k)}).
        \end{align*}
    
\item     The set $\mathcal{S}$ in Lemma \ref{LEMMAEXISTENCE} is in fact the $c$-superdifferential of the $c$-concave potential $\phi$ and is called a \textit{minimizing set}  which contains the support of all optimal plans for \eqref{TAGMKP}.
        \end{enumerate}
\end{remark}

Before proceeding further, we recall a newly established  tool that could assist in demonstrating the uniqueness of the optimal plan of \eqref{TAGMKP}. Indeed, in \cite{MOMENI-RIFFORD} a criterion regarding the uniqueness of the optimal plan of \eqref{TAGMKP} is obtained which is based on the notion of $c$-extreme minimizing sets. In this regard, we introduce the following set-valued functions,
\begin{align}
    &F:X\to 2^Y,\quad F(x)=\left\{y\in Y\ :\ (x,y)\in \mathcal{S} \right\},\label{F-function}\\
    &f:X\to 2^Y,\quad f(x)=\argmax\{c(x,y)\ :\  y\in F(x)\},\nonumber
\end{align}
where the set $\mathcal{S}$ is given by Lemma \ref{LEMMAEXISTENCE}. Domains of these functions are naturally defined as follows
\begin{align*}
    \text{Dom}(F)=\left\{x\ :\ F(x)\neq\emptyset\right\},\quad     \text{Dom}(f)=\left\{x\ :\ f(x)\neq\emptyset\right\}.
\end{align*}
Now we are prepared to give the definition of a $c$-extreme minimizing set.
\begin{definition}\label{DEF-C-EXT}
We say the minimizing set $ \mathcal{S}$ of the problem  \eqref{TAGMKP} is $c$-extreme if there exist $\mu$- and $\nu$-full measures subsets $M$ and $N$ of $X$ and $Y$, respectively, such that
\begin{enumerate}
    \item $\text{Dom}(F)\cap M = \text{Dom}(f)\cap M$,
    \item For all distinct $x_1$ and $x_2$ in $\text{Dom}(F)\cap M$ and for all $y_i\in f(x_i)$, for $ i=1,2$, we have that
    \begin{align*}
        \left(F(x_1)\setminus\{y_1\}\right)\cap \left(F(x_2)\setminus\{y_2\}\right)\cap N=\emptyset .
    \end{align*}
\end{enumerate}
\end{definition}
The following theorem states the sufficient condition for \eqref{TAGMKP} to have a unique optimal plan, which is concentrated on the disjoint union of the graph and anti-graph of two measurable maps. Interested readers may refer to \cite[Theorem 2.10]{MOMENI-RIFFORD} and \cite[Corollary 2.3]{MOMENI-DOUBLY} for relevant findings and proof techniques.

\begin{theorem}\label{TH-C-EXT}
Let $X$ and $Y$ be smooth closed manifolds equipped with the Borel probability measures $\mu$ and $\nu$, respectively. Assume that $c : X\times Y \to\mathbb{R}$ is a continuous function and $\mathcal{S} \subseteq X \times Y$ is a minimizing set. If $\mathcal{S}$ is $c$-extreme, then there exists a unique $\lambda\in\Pi(\mu,\nu)$ with $\lambda(\mathcal{S})=1$. Moreover, the measure $\lambda$ is the extreme point of the set $\Pi(\mu,\nu)$.   
\end{theorem}
A more general statement of Definition \ref{DEF-C-EXT} and Theorem \ref{TH-C-EXT} can be expressed in the case that the target space $Y$ is partitioned. To state them, we recall that a family $P=\{Y_k\}_{k=1}^K$ is a \textit{Borel ordered partition} of $Y$ if each $Y_k$ is a Borel set, they are pairwise disjoint, and we have $Y=\bigcup_{k=1}^K Y_k$. So, let us consider the problem \eqref{TAGMKP} and its dual \eqref{TAGDMKP} for such a target space $Y$, and define the following set-valued functions
\begin{align*}
&\kappa :\text{Dom}(F)\to\{1,\ldots,K\},\quad \kappa(x)=\min\{k\in\{1,\ldots,K\}\ : F(x)\cap Y_k\neq\emptyset\},\\
&f_P:X\to 2^Y,\quad f_P(x)=\arg\max\{c(x,y)\ :\  y\in F(x)\cap Y_{\kappa(x)}\}.
\end{align*}
\begin{remark}\label{REM-KAPPA}
For $Y$ with a Borel ordered partition $P=\{Y_{k}\}_{k=1}^K$, it is worth noting that according to the definitions of $\kappa$ and $f_P$, if  $y\in Y_i$ and $y^*\in Y_j$ satisfy
\begin{align*}
   y \in f_P(x)\quad\text{and}\quad  y^* \in F(x),\quad \text{for some $x\in X$},
\end{align*}
then we have $i\leq j$.
\end{remark}

The next definition generalizes the idea of a $c$-extreme minimizing set to cases where the target space is divided into several parts. This expansion helps us understand optimization problems better in situations where the target space is split into different regions. By looking at these partitions, we can analyze how optimization works in each region separately, giving us a clearer picture of the problem overall.
\begin{definition}\label{DEF-CP-EXT}
Let $P=\{Y_k\}_{k=1}^K$ be a Borel ordered partition of $Y$. Then the minimizing set $\mathcal{S}$ of the problem \eqref{TAGMKP} is called $(c,P)$-extreme if there are $\mu$- and $\nu$-full measure subsets $M$ and $N$ of $X$ and $Y$, respectively, such that
\begin{enumerate}
    \item $\text{Dom}(F)\cap M = \text{Dom}(f_P)\cap M$,
    \item For all distinct $x_1$ and $x_2$ in $\text{Dom}(F)\cap M$ and for all $y_i\in f_P(x_i), i=1,2$, we have that
    \begin{align*}
        \left(F(x_1)\setminus\{y_1\}\right)\cap \left(F(x_2)\setminus\{y_2\}\right)\cap N=\emptyset .
    \end{align*}
\end{enumerate}
\end{definition}

Now, we recall  a practical method for demonstrating the uniqueness of the optimal solution in \eqref{TAGMKP} when the target space $Y$ consists of disjointed Borel sets (see \cite[Theorem 2.12]{MOMENI-RIFFORD}). In fact, the upcoming theorem generalizes the result of Theorem \ref{TH-C-EXT} in the case of a partitioned space $Y$.

\begin{theorem}\label{TH-CP-EXT}
Let $X$ and $Y$ be smooth closed manifolds equipped with probability measures $\mu$ and $\nu$, respectively. Assume that $c : X\times Y \to\mathbb{R}$ is a continuous function and $\mathcal{S} \subseteq X \times Y$ is a minimizing set. Moreover, let $P=\{Y_k\}_{k=1}^K$ be a Borel ordered partition of $Y$. If $\mathcal{S}$ is $(c,P)$-extreme, then there exists a unique $\lambda\in\Pi(\mu,\nu)$ with $\lambda(\mathcal{S})=1$. Moreover, the measure $\lambda$ is the extreme point of the set $\Pi(\mu,\nu)$.

\end{theorem}
In \cite{MOMENI-CHAR}, a criterion is presented for recognizing the support of optimal plans of Monge-Kantorovich problems. This criterion is based on the notion of an $m$-twist condition, which guarantees that the support of optimal plans is concentrated on the graph of several maps. Here, we briefly describe this notion and its consequences. First, we need to clarify of what we mean by concentrating on the graphs of several maps. 

\begin{definition}\label{DEF-M-G}
    A measure $\lambda\in \Pi(\mu,\nu)$ is concentrated on the union of the graph of $m$ measurable functions $T_i:X\to Y$, for $i=1,\ldots,m$, if there exists $m$ measurable maps $\{\alpha_i\}_{i=1}^m$ such that
    \begin{align*}
        \alpha_i :X\to [0,1],\quad \sum_{i=1}^m \alpha_i(x)=1,\quad\text{for $\mu$-a.e.}\ x\in X,
    \end{align*}
    for which we have
    \begin{align}\label{M-GH-LAM}
        \lambda(A\times B)=\sum_{i=1}^m\int_A \alpha_i(x)\chi_{B}(T_i(x))\, d\mu(x),\quad A\times B\in \mathcal{B}_{X}\times \mathcal{B}_{Y}.
    \end{align}
\end{definition}

An additional definition is required to adequately prepare for the subsequent result.

\begin{definition}\label{DEF-M-TWIST}
For $m\in\mathbb{N}$, we say that the cost function $c:X\times Y\to\mathbb{R}$ satisfies the $m$-twist condition, if the following set has at most $m$ elements
\begin{align*}
    L(x_0,y_0)=\left\{y\in Y\ :\ \nabla_xc(x_0,y)=\nabla_xc(x_0,y_0) \right\},
\end{align*}
 for $\mu\otimes\nu$-almost every $(x_0,y_0)\in X\times Y$.
\end{definition}

The following theorem determines the circumstances under which the optimal plans of the problem \eqref{TAGMKP} are concentrated on the union of the graph of $m$ measurable functions, provided that the first marginal has no atom.

\begin{theorem}\label{M-TWIST-TH}
Let $X$ be a complete separable Riemannian manifold and $Y$ be a Polish space equipped with Borel probability measures $\mu$ and $\nu$, respectively. Let $c : X\times Y \to \mathbb{R}$ be a bounded continuous cost function such that
\begin{enumerate}
    \item the cost function $c$ satisfies the $m$-twist condition,
    \item the marginal $\mu$ is non-atomic, and any $c$-concave function on $X$ is differentiable $\mu$-almost everywhere on its domain.
\end{enumerate}
Then each optimal plan $\lambda$ of \eqref{TAGMKP} is concentrated on the union of the graph of $m$ measurable functions, that is, $\lambda$ is of the form \eqref{M-GH-LAM}.
\end{theorem}

\begin{remark}
The notion of $c$-extreme and $(c,P)$-extreme minimizing sets can be easily extended to the multi-marginal case. For instance, consider $(X,\mu)$, $(Y,\nu)$, and $(Z,\gamma)$ as probability spaces, with the cost function defined as $c:X\times Y\times Z\to\mathbb{R}$. We can then treat it as a function $c:W\times Z\to\mathbb{R}$, where $W=X\times Y$. In cases where we can apply Theorem \ref{TH-CP-EXT}, the optimal plan $\lambda\in \Pi(\mu,\nu,\gamma)$ becomes an extreme point of $\Pi(\lambda^{XY},\gamma)$, where $\lambda^{XY}$ denotes the restriction of $\lambda$ to $X\times Y$. Furthermore, Definitions \ref{DEF-M-TWIST} and \ref{DEF-M-G}, along with the result of Theorem \ref{M-TWIST-TH}, can be extended to the cost function $c$ using the same approach.
\end{remark}

\section{Strictly convex costs on multi-layers target spaces}\label{SECTIONMULTI-LAYER-STRCONV}

In this section, we study cost functions of the form $c(x,y)=h(x-y)$, where $\displaystyle{h:\mathbb{R}^{n+1}\to\mathbb{R}^+}$ is a strictly convex function. The problem \eqref{TAGMKP} with such a cost function has been studied in \cite{AbdellaouiHeinich,BRENIER,CAFFARELLI,FIGALLI,GANGBOMCCANN,GANGBOMCCANN2,KNOTT,RUSCH,SMITH}. It has been established that when the first marginal is absolutely continuous with respect to the Lebesgue measure on $\mathbb{R}^{n+1}$, there is a unique measurable map $T:X\to Y$ that induces the unique optimal plan. It is shown that the geometry of the domain and the role of $c$-concave functions are undeniable.


In a seminal paper, Brenier proved that for the squared Euclidean distance $c(x, y) = |x - y|^2$, there exists a unique optimal map, identified as the gradient of a convex function. Utilizing the Kantorovich dual problem, one can handle both non-uniform mass distributions in $\mathbb{R}^{n+1}$ and uniform distributions on sets $X$ and $Y$, under the assumption that the total masses are equal. This result, reminiscent of Riemann's mapping theorem, leads to consequences such as a polar factorization theorem for vector fields and a Brunn-Minkowski inequality for measures. The implications of these discoveries raise fundamental questions about the features of the cost function that determine the existence and uniqueness of optimal maps, the geometric properties characterizing maps for other costs, and the potential fruitful applications of their geometry \cite{AbdellaouiHeinich,BRENIER,CAFFARELLI,KNOTT,RUSCH,SMITH}.

Furthermore, in \cite{GANGBOMCCANN}, the authors dealt with two important classes of cost functions. Namely, $c(x,y)=h(x-y)$ with $h$ strictly convex, and $c(x,y)=l(|x-y|)$ with $l\geq 0$ strictly concave. For convex cost functions, a theory parallel to that for the squared Euclidean distance has been developed: the optimal map exists and is uniquely characterized by its geometry. This map depends explicitly on the gradient of the cost, or rather on its inverse map $(\nabla h)^{-1}$, which indicates why strict convexity or concavity should be essential for uniqueness. Although explicit solutions are more awkward to obtain, we have no reason to believe that they should be any worse behaved than those for squared Euclidean distance (see e.g., the regularity theory developed in \cite{CAFFARELLI,FIGALLI}). Moreover, in \cite{GANGBOMCCANN2}, the authors considered the quadratic problem for the case where the marginals are defined on the boundary of bounded convex domains $X$ and $Y$ in $\mathbb{R}^{n+1}$. It has been shown that when the first marginal is absolutely continuous with respect to the $n$-dimensional Hausdorff measure on $\partial X$, then the optimal plan is unique and more regularity results are obtained under extra conditions such as uniform convexity of the domains and absolute continuity of both marginals with special densities.

To continue working on such problems, in this section, we address the two-marginal Monge-Kantorovich optimization problem
\begin{align}\label{STR-CONV-NORMAL}
    \inf\left\{\int_{X\times Y}h(x-y)\,d\lambda\ :\ \lambda\in\Pi(\mu,\nu)\right\},
\end{align}
where the function $h:\mathbb{R}^{n+1}\to\mathbb{R}^+$ is strictly convex, and the domains $X$ and $Y$ are special $n$-dimensional differentiable manifolds of $\mathbb{R}^{n+1}$. Let us describe the domains $X$ and $Y$. Assume that $K\in\mathbb{N}$, and $\overline{X},\overline{Y}_{1},...,\overline{Y}_{K} \subseteq\mathbb{R}^{n}$ are compact subsets of $\mathbb{R}^{n}$. Let $\overline{x}\in\mathbb{R}$, and $\{\overline{y}_1,..., \overline{y}_K\}\subseteq \mathbb{R}$ be a set of distinct real numbers. Set
\begin{align*}
    X=\overline{X}\times\{\overline{x}\},\quad\text{and}\quad Y=\bigcup_{k=1}^K \left(\overline{Y}_{k}\times \{\overline{y}_k\}\right).
\end{align*}
Therefore, the subsets $X,Y\subseteq\mathbb{R}^{n+1}$ are $n$-dimensional sub-manifolds.
Our goal is to show that due to the geometry of the domains, the solution to problem \eqref{STR-CONV-NORMAL} is unique, but it is not induced by a single map. More precisely, we have the following theorem.

\begin{theorem}\label{TH31}
    Assume that $\overline{X},\overline{Y}_k\subseteq \mathbb{R}^n$ are compact subsets, for $k=1,\ldots,K$ with $K\in\mathbb{N}$. Put
    \begin{align*}
        X=\overline{X}\times \{\overline{x}\},\quad Y=\bigcup_{k=1}^K\left(\overline{Y}_k\times\{\overline{y}_k\}\right),
    \end{align*}
where $\overline{x}\in\mathbb{R}$, and $\{\overline{y}_1,\ldots,\overline{y}_K\}$ is a set of distinct elements in $\mathbb{R}$ .
    Equip $X$ and $Y$ with the Borel probability measures $\mu$ and $\nu$, respectively as follows
    \begin{itemize}
        \item $\mu\ll\mathcal{L}^n$.
        \item If $K\geq 2$, then $\nu=\displaystyle{\sum_{k=1}^K t_k\nu_k}$, where $\nu_k$'s are Borel probability measures on $\overline{Y}_k\times \{\overline{y}_k\}$, for $k=1,\ldots,K$,  with $\nu_k\ll\mathcal{L}^n$, for $k=2,\ldots,K$, and $\displaystyle{\sum_{k=1}^K t_k=1}$. 
    \end{itemize}
    Let  $h:\mathbb{R}^{n+1}\to\mathbb{R}^+$ be a differentiable  strictly convex function. Then, the optimization problem \eqref{STR-CONV-NORMAL} admits a unique solution, which is concentrated on the graph of $K$ functions $T_k:X\to Y$, for $k=1,\ldots,K$.
\begin{proof}
    Existence of the minimizer of \eqref{STR-CONV-NORMAL} is guaranteed because of the compactness of $X$ and $Y$.
        To show the uniqueness, first recall that
        \begin{align*}
            \mathcal{S}=\left\{(x,y)\ :\ h(x-y)=\phi(x)+\psi(y) \right\},
        \end{align*}
        where $\phi$ and $\psi$ are solutions (potentials) to the dual problem of \eqref{STR-CONV-NORMAL}. We shall prove that the minimizing set $\mathcal{S}$ of \eqref{STR-CONV-NORMAL} is $(c,P)$-extreme, where $P=\{Y_k\}_{k=1}^K$ is a Borel ordered partition for $Y$. Here, we have set 
        \begin{align*}
            Y_k=\overline{Y}_k\times\{\overline{y}_k\},\quad k=1,\ldots,K.
        \end{align*}
        Consider the  set-valued function
        \begin{align*}
            &F:X\to 2^Y,\quad F(x)=\{y\in Y\ :\ (x,y)\in \mathcal{S}\}.
        \end{align*}  
        Moreover, recall that
        \begin{align*}
            \kappa : \text{Dom}(F)\to\{1,\ldots,K\},\quad \kappa(x)=\min\{k\in\{1,\ldots,K\}\ :\ F(x)\cap Y_k\neq \emptyset \},
        \end{align*}
        and define
        \begin{align*}
            &f_P:X\to 2^Y.\quad f_P(x)=\argmax\{h(x-y)\ : y\in F(x)\cap Y_{\kappa(x)}\}.
        \end{align*}
        We must show that the following criteria are satisfied
\begin{enumerate}
    \item $\text{Dom}(F)\cap N=\text{Dom}(f)\cap N$,
    \item $F(x_1)\backslash \{y_1\}\cap F(x_2)\backslash \{y_2\}\cap M=\emptyset,\quad\forall y_i\in f_P(x_i),\ \forall x_i\in M$ with $x_1\neq x_2$,
\end{enumerate}
where the $\mu-$ and $\nu-$ full measure subsets $M\subseteq X$ and $N\subseteq Y$ should be determined.
We note that the potentials are locally Lipschitz; therefore, because of the fact that $\mu\ll\mathcal{L}^n$ and $\nu_k\ll\mathcal{L}^n$, for $k=2,\ldots,K$, the following subsets are of $\mu$- and $\nu$-full measures
\begin{align*}
    M=\{x\in X\ : \ \nabla\phi(x)\ \text{exists}\},\quad N=Y_1\cup\left(\bigcup_{k=2}^K\left\{y\in Y_k\ :\ \nabla\psi(y)\ \text{exists}\right\}\right). 
\end{align*}
The first condition is satisfied by the compactness of $Y$ and continuity of $h$. To show the second condition, we assume that the intersection is non-empty, and we will come up with a contradiction. So, let $y^*$ be in the intersection. Note that in this case we have 
\begin{align*}
    (x_1,y_1),\;(x_2,y_2),\; (x_1,y^*),\;(x_2,y^*)\in\mathcal{S}.
\end{align*}
We consider the two following cases.

\begin{enumerate}
    \item[Case 1:] $y^*\in Y_1= \overline{Y}_1\times\{\overline{y}_1\}$. In this case, according to Remark \ref{REM-KAPPA}, we must have that $y_1,y_2\in Y_1$. Since $x_1\in M$, then by the differentiability of $\phi$ at $x_1$, the fact that $(x_1,y_1),\;(x_1,y^*)\in\mathcal{S}$, and with the help of Lemma \ref{LEM-CRITICAL}, we obtain that
\begin{align}\label{Eq21}
e_{n+1}=\alpha(\nabla h(x_1-y_1)-\nabla h(x_1-y^*)),\quad\text{for some}\;\alpha\in\mathbb{R}\backslash\{0\},
\end{align}
where $e_{n+1}$ denotes the $(n+1)$-th standard coordinate vector. 
Now, by multiplying both sides of \eqref{Eq21} by $y^*-y_1$, we get
\begin{align}\label{E24}
    0=\langle e_{n+1},y^*-y_1\rangle=\alpha\langle \nabla h(x_1-y_1)-\nabla h(x_1-y^*), y^*-y_1\rangle\neq 0,
\end{align}
which is a contradiction.
 \item[Case 2:] $y^*\in\bigcup_{k=2}^K \left(\overline{Y}_k\times\{\overline{y}_k\}\right)\cap N$. In this case, $\nabla\psi(y^*)$ exists, and $(x_1,y^*),\;(x_2,y^*)\in\mathcal{S}$. Consequently, with the same approach we obtain that
\begin{align}\label{Eq22}
e_{n+1}=\alpha(\nabla h(x_1-y^*)-\nabla h(x_2-y^*)),\quad\text{for some}\;\alpha\in\mathbb{R}\backslash\{0\}.
\end{align}
Multiplying both sides of \eqref{Eq22} by $x_1-x_2$ gives us the same contradiction as \eqref{E24}.
\end{enumerate}      

Hence, the second condition is satisfied, and the minimizing set $\mathcal{S}$ is $(c,P)$-extreme. Therefore, the optimization problem \eqref{STR-CONV-NORMAL} admits a unique solution by Theorem \ref{TH-CP-EXT}. To show that the support of the optimal plan is concentrated on the union of the graph of $K$ functions $T_k:X\to Y$, for $k=1,\ldots,K$, we will demonstrate that the cost function $c(x,y)=h(x-y)$ satisfies the $K$-twist condition. We shall prove that for any $(x_0,y_0)\in M\times N$, the following set has at most $K$ elements
\begin{align*}
    L(x_0,y_0)=\left\{y\in F(x_0) \ :\ \nabla h(x_0-y)=\nabla h(x_0-y_0) \right\}.
\end{align*}
We already know that $y_0\in L(x_0,y_0)$. Suppose otherwise, and let $L(x_0,y_0)$ contain at least $K+1$ distinct elements $\{y_i\}_{i=0}^K$ of $F(x_0)$. Then there exists $k_0\in\{1,\ldots, K\}$ such that $\overline{Y}_{k_0}\times \{\overline{y}_{k_0}\}$ contains at least two elements $y_{k_1}$ and $y_{k_2}$ from $\{y_i\}_{i=0}^K$. Now, due to the fact that $(x_0,y_{k_1}),\; (x_0,y_{k_2})\in \mathcal{S}$ and the differentiability of $\phi$ at $x_0$, we have
\begin{align}\label{Eq23}
    e_{n+1}=\alpha(\nabla h(x_0-y_{k_1})-\nabla h(x_0-y_{k_2})),\quad\text{for some}\;\alpha\in\mathbb{R}\backslash\{0\}.
\end{align}
Multiplying both sides of \eqref{Eq23} by $y_{k_2}-y_{k_1}$ leads to the same contradiction as \eqref{E24}. Thus, via Theorem \ref{M-TWIST-TH}, we obtain the result.

    \end{proof}
\end{theorem}

It is worth noting that in the case where $K=1$, the problem reduces to an optimization problem on $\mathbb{R}^n$, having a unique solution concentrated on the graph of a single optimal transport map \cite{GANGBOMCCANN}.

\subsection{Optimality of Theorem \ref{TH31} for the quadratic cost}

In this subsection, we refine Theorem \ref{TH31} for the case where $h(x-y)$ is the quadratic cost function. Moreover, by providing two examples, we justify the optimality of our assumptions. In fact, for the quadratic function $h(x-y)=|x-y|^2$, we can relax the geometry of our domains $X$ and $Y$ as stated in Theorem \ref{TH31}.

\begin{theorem}\label{TH32}
    Assume that $P$ and $P_k$, for $k=1,\ldots,K$ with $K\in\mathbb{N}$, are planes in $\mathbb{R}^{n+1}$ with normal vectors $\mathrm{n}$ and $\mathrm{n}_k$, and passing through the points $\overline{x}_0\in \mathbb{R}^{n+1}$ and $\overline{y}_k\in \mathbb{R}^{n+1}$, respectively.  Moreover, assume that
    \begin{align}\label{PERP}
        \langle\mathrm{n}, \mathrm{n}_k\rangle\neq 0,\quad k=1,\ldots,K.
    \end{align}
    Let $X, Y_k\subseteq \mathbb{R}^{n+1}$ be compact subsets of the planes $P$ and $P_k$, respectively, and set $Y=\bigcup_{k=1}^K Y_k$. Equip $X$ and $Y$ with the Borel probability measures $\mu$ and $\nu$, respectively, as follows
    \begin{itemize}
        \item $\mu\ll\mathcal{L}^n$.
        \item  If $K\geq 2$, then $\displaystyle{\nu=\sum_{k=1}^K t_k\nu_k}$ where $\nu_k$'s are Borel probability measures on $Y_k$ with $\nu_k\ll\mathcal{L}^n$, for $k=2,\ldots,K$, and $\displaystyle{\sum_{k=1}^K t_k=1}$. 
    \end{itemize}
Then, the optimization problem 
   \begin{align}\label{2-MARG-QUAD-PROB}
    \inf\left\{\int_{X\times Y}|x-y|^2\,d\lambda\ :\ \lambda\in\Pi(\mu,\nu)\right\},
\end{align}
admits a unique solution which is concentrated on the graph of $K$ functions $T_k:X\to Y$, for $k=1,\ldots,K$.
    \begin{proof}
The proof is an adaptation of the proof of Theorem \ref{TH31}. In fact, consider the set-valued maps $F$, $\kappa$, and $f_P$ for the case of quadratic cost and the domains $X$ and $Y$ with the property \eqref{PERP} and a partition $P=\{Y_k\}_{k=1}^K$ for $Y$. Define
\begin{align*}
    M=\{x\in X\ : \ \nabla\phi(x)\ \text{exists}\},\quad N=Y_1\cup\left(\bigcup_{i=2}^K\left\{y\in Y_i\ :\ \nabla\psi(y)\ \text{exists}\right\}\right). 
\end{align*}
where $\phi$ and $\psi$ are the potentials as before (existence of such continuous potentials are always guaranteed in this framework).  We just need to modify our computations in establishing the second condition in showing that the minimizing set $\mathcal{S}$ is $(c,P)$-extreme. In fact, assume that 
\begin{align*}
    y^*\in F(x_1)\setminus\{y_1\}\cap F(x_2)\setminus\{y_2\},
\end{align*}
for some $x_1\neq x_2\in M$ and $y_i\in f_P(x_i)$, for $i=1,2$. Note that in this case we have 
\begin{align*}
    (x_1,y_1),\;(x_2,y_2),\; (x_1,y^*),\;(x_2,y^*)\in\mathcal{S}.
\end{align*}
Similarly, we consider the two following cases.
        \begin{enumerate}
            \item[Case 1:] $y^*\in Y_1$. Through Remark \ref{REM-KAPPA} we must have that $y_1,y_2\in Y_1$. Since $x_1\in M$, again via Lemma \ref{LEM-DIFF} we obtain that
\begin{align}\label{Eq26}
\mathrm{n}=\alpha(y^*-y_1),\quad\text{for some}\;\alpha\in\mathbb{R}\backslash\{0\}.
\end{align}
Now, by multiplying both sides of \eqref{Eq26} by $\mathrm{n}_1$, we have that
\begin{align*}
    0\neq \langle \mathrm{n}, \mathrm{n}_1\rangle=\alpha\langle y^*-y_1, \mathrm{n}_1\rangle= 0,
\end{align*}
which  is a contradiction.
 \item[Case 2:] $y^*\in\bigcup_{i=2}^K \left(Y_i\cap N\right)$. In this case $\nabla\psi(y^*)$ exists. Together with the fact that $(x_1,y^*),\;(x_2,y^*)\in\mathcal{S}$, we have
\begin{align}\label{Eq27}
\mathrm{n}_{k_0}=\alpha(x_1-x_2),\quad\text{for some}\;\alpha\in\mathbb{R}\backslash\{0\},
\end{align}
where $k_0$ is the index of the $Y_k$ that contains $y^*$. Again, by multiplying both sides of \eqref{Eq27} by $\mathrm{n}$, the normal vector of the plane $P$ containing $X$, we obtain that
\begin{align}
    0\neq \langle \mathrm{n}_{k_0},\mathrm{n}  \rangle=\alpha\langle x_1-x_2, \mathrm{n} \rangle= 0,
\end{align}
which is a contradiction.
\end{enumerate}
Therefore, the minimizing set $\mathcal{S}$ is $(c,P)$-extreme. Thus, the solution to \eqref{2-MARG-QUAD-PROB} is unique. The fact that the the support of the unique optimal plan is concentrated on the union of the graph of functions $T_k:X\to Y$, for $k=1,\ldots,K$, can be obtained similarly. 
    \end{proof}
    \end{theorem}
 We remark that for the case that $K=1$ our result in the previous Theorem  coincides with the result in \cite{GANGBOMCCANN2}.
In the upcoming examples, we aim to illustrate potential cases where the uniqueness of optimal plans could fail, particularly when the sets $X$ and $Y$ exhibit orthogonal characteristics or when the second marginal is not absolutely continuous on the second layer $Y_2$ (and the layers afterward) of the target space $Y$.
\begin{example}
In this example, we demonstrate that by removing the requirement of absolute continuity on $Y$, the uniqueness of the solution may no longer hold.
Let 
\begin{align}
X=[0,1]\times \{0\} \subset \R^2,\quad \text{and}\quad Y= Y_1\cup Y_2  \subset \R^2\quad \text{where} \quad Y_i=[0,1]\times \{(-1)^{i+1}\}.
\end{align}
Define the marginal measures $\mu$ and $\nu$ on $\R^2$ as follows,
\[\int_{\R^2} f(x_1,x_2)\,d\mu=\int_0^1 f(x_1,0) \,d x_1,  \]
and
\[\int_{\R^2} f(y_1,y_2)\,d\nu=\frac{1}{2}f(1,1)+\frac{1}{2}f(1,-1),\]
for all continuous and bounded functions $f.$
Note that $\mu$ is supported on $X$, and is absolutely continuous with respect to the one dimensional Lebesgue measure.  Note also that  $\nu$ is supported on $Y$,  and its restriction to either $Y_1$ or $Y_2$ is not absolutely continuous   with respect to  the one dimensional Lebesgue measure.
For the cost function $c: X \times Y\to \R$ defined by $c(x,y)=|x-y|^2,$ assume that $\lambda_0 \in \Pi(\mu,\nu)$ minimizes 
\[\inf\left\{\int_{X\times Y} |x-y|^2 d\, \lambda\ :\ \lambda \in \Pi(\mu,\nu)\right\}.\]
Note that for each $(x,y)\in \text{Spt}(\lambda_0)$ with $x=(x_1,x_2) \in \text{Spt}(\mu)$
and $y=(y_1, y_2)\in \text{Spt}(\nu)$ we have that 
\[|x-y|^2=|x_1-1|^2+1.\]
It then follows that
\[\int |x-y|^2 d\, \lambda_0=\int (|x_1-1|^2+1) d\lambda_0=\int_0^1 (|x_1-1|^2+1) dx_1.\]
We shall now construct two different optimal plans.  Define the maps $T_1: X \to Y$ and $T_2:X \to Y$ as follows,
\begin{equation}\label{e1}
    T_1(x_1, x_2)= \left\{\begin{aligned}
    & (1,1), \qquad x_1 \in [0, \frac{1}{4})\cup [\frac{1}{2}, \frac{3}{4})\\
     &\\
    & (1,-1), \qquad x_1 \in [\frac{1}{4}, \frac{1}{2})\cup [\frac{3}{4},1]
     \end{aligned}
     \right.,
	\end{equation}
 and 

 \begin{equation}\label{e2}
    T_2(x_1, x_2)= \left\{\begin{aligned}
     &(1,-1), \qquad x_1 \in [0, \frac{1}{4})\cup [\frac{1}{2}, \frac{3}{4})\\
    & \\
     &(1,1), \qquad x_1 \in [\frac{1}{4}, \frac{1}{2})\cup [\frac{3}{4},1]
     \end{aligned}
     \right. .
	\end{equation}
It can be easily deduced that $T_i\#\mu=\nu$  for $i=1,2.$  Moreover, for $i=1,2$ we have that 

\[\int_X |x-T_i(x)|^2 d\mu=\int_0^1 (|x_1-1|^2+1) dx_1=\int |x-y|^2 d\, \lambda_0.\]

This shows that the two distinct  plans $\lambda_1=(\text{id}\times T_1)\#\mu$ and $\lambda_2=(\text{id}\times T_2)\#\mu$ are indeed optimal. Therefore, we do  not have uniqueness in the case where  we relax the absolute continuity assumption on $\nu$ required in Theorem \ref{TH31}.
\end{example}

\begin{example}
In this example we show that if $X $ and $Y$ are perpendicular  then the solution may not be unique.  Let $X=[0,1]\times \{0\}$ and $Y=\{0\}\times [0,1]$ with probability measure $\mu$ supported on $X$ and probability measure $\nu$ supported on $Y$.  Then for the cost function $c(x,y)=|x-y|^2$,  and any $\lambda\in\Pi(\mu, \nu)$ we have that 
\[\int_{X\times Y} |x-y|^2 \, d\lambda=\int_{X\times Y} (|x_1|^2+|y_2|^2) \, d\lambda=\int_X |x_1|^2 \, d\lambda+\int_Y |y_2|^2 \, d\gamma=\int_X |x_1|^2 \, d\mu+\int_Y |y_2|^2 \, d\nu. \]

This shows that any plan $\gamma \in \Pi(\mu, \nu)$ is indeed an optimal plan,  and therefore the uniqueness fails. 
\end{example}

\section{Multi-marginal quadratic cost on multi-layers target spaces}\label{SECTIONMULTI-LAYER-QUAD23}
This section is devoted to studying the minimization problem
\begin{align*}
    \inf\left\{\int\left(|x-y|^2+|x-z|^2+|y-z|^2\right)\,d\lambda \ : \ \lambda\in\Pi(\mu,\nu,\gamma)\right\},
\end{align*}
where the domains $X$, $Y$, and $Z$ are $n$-dimensional sub-manifolds of $\mathbb{R}^{n+1}$, equipped with the Borel probability measures $\mu$, $\nu$, and $\gamma$, respectively. For simplicity in computations, we consider the following equivalent maximization problem instead
\begin{align}\label{3-MAR-QUAD-LAYOUT}
        \sup\left\{\int\left(\langle x,y\rangle+\langle x,z\rangle+\langle y,z\rangle\right)\,d\lambda\ : \ \lambda\in\Pi(\mu,\nu,\gamma)\right\}.
\end{align}
Our goal is to establish the uniqueness of the solution to \eqref{3-MAR-QUAD-LAYOUT} by utilizing the theory of $(c,P)$-extreme maximizing sets and a reduction argument. Here, we present the main result of this section.
\begin{theorem}\label{TH-3MARG-QUAD-MUTI-LAY}
Let $\overline{X}$, $\overline{Y}_k$, and $\overline{Z}_l$ be compact subsets of $\mathbb{R}^n$, and let $\overline{x}\in\mathbb{R}$, $\{\overline{y}_1,\ldots,\overline{y}_K\}$, and $\{\overline{z}_1,\ldots,\overline{z}_L\}$ be distinct elements of $\mathbb{R}$, with $K,L\in\mathbb{N}$. Define
\begin{align*}
    X=\overline{X}\times\{\overline{x}\},\quad Y=\bigcup_{k=1}^K\left(\overline{Y}_k\times\{\overline{y}_k\}\right),\quad\text{and}\quad Z=\bigcup_{l=1}^L\left(\overline{Z}_l\times\{\overline{z}_l\}\right).
\end{align*}
Equip $X$, $Y$, and $Z$ with the Borel probability measures $\mu$, $\nu$, and $\gamma$, respectively, as follows
\begin{itemize}
    \item $\mu\ll\mathcal{L}^n$,
    \item If $K\geq 2$, then $\displaystyle{\nu=\sum_{k=1}^K t_k\nu_k}$, where $\nu_k$'s are Borel probability measures on $\overline{Y}_k\times\{\overline{y}_k\}$ such that $\nu_k\ll\mathcal{L}^n$ for $k=2,\ldots,K$, and $\displaystyle{\sum_{k=1}^K t_k=1}$.
    \item If $L\geq 2$, then ${\displaystyle\gamma=\sum_{l=1}^L s_l\gamma_l}$, where $\gamma_l$'s are Borel probability measures on $\overline{Z}_l\times\{\overline{z}_l\}$ such that $\gamma_l\ll\mathcal{L}^n$ for $l=2,\ldots,L$, and $\displaystyle{\sum_{l=1}^L s_l=1}$.
\end{itemize}
Then, the maximization problem \eqref{3-MAR-QUAD-LAYOUT} admits a unique solution, concentrated on the graphs of $K\times L$ maps $T_k:X\to Y\times Z$.

    \begin{proof}
Recall that the dual problem to \eqref{3-MAR-QUAD-LAYOUT}  is of the following form,
        \begin{align}\label{3-MAR-QUAD-LAYOUT-DUAL}
            \min\left\{\int_X\phi_1(x)\,d\mu(x)+\int_Y\phi_2(y)\,d\nu(y)+\int_Z\phi_3(z)\,d\gamma(z)\ :\ \phi_1(x)+\phi_2(y)+\phi_3(z)\geq c(x,y,z) \right\},
        \end{align}
       where the supremum is taken over all triples $(\phi_1,\phi_2,\phi_3)\in L_1(X,\mu)\times L_1(Y,\nu)\times L_1(Z,\gamma)$, that their summation at each point $(x,y,z)\in X\times Y\times Z$ is not dominated by the cost $c(x,y,z)$. By Theorem \ref{TH-DUALATTAIN}, existence of the solutions to the problem \eqref{3-MAR-QUAD-LAYOUT} and \eqref{3-MAR-QUAD-LAYOUT-DUAL} is guaranteed. Let us work with the potentials $\phi_1$, $\phi_2$, and $\phi_3$, and the following  maximizing set $\mathcal{S}$ associated with the problem \eqref{3-MAR-QUAD-LAYOUT},
\begin{align*}
    \mathcal{S}=\left\{(x,y,z)\ : \ \phi_1(x)+\phi_2(y)+\phi_3(z)=c(x,y,z) \right\}.
\end{align*}
Using the potential $\phi_3$, we associate the following 2-marginal problem to the 3-marginal problem \eqref{3-MAR-QUAD-LAYOUT},
        \begin{align}\label{ASSOC-2MARG-TO-3-MAR-QUAD-LAYOUT}
        \sup\left\{\int\left(\langle x,y\rangle+\xi(x+y)\right)\,d\tau\ : \ \tau\in\Pi(\mu,\nu)\right\},
\end{align}
where
\begin{align*}
    \xi(x+y)=\sup\left\{\langle x+y,z\rangle-\phi_3(z)\ : z\in Z\right\}.
\end{align*}
Let $\lambda\in\Pi(\mu,\nu,\gamma)$ be an optimal plan for  \eqref{3-MAR-QUAD-LAYOUT}. By Proposition \ref{PROP-OPTIMAL}, we know that $\lambda_{XY}$, the restriction of $\lambda$ on $X\times Y$ is an optimal plan of \eqref{ASSOC-2MARG-TO-3-MAR-QUAD-LAYOUT} and via Lemma \ref{LEM-MINI-SSJ}, the set $\pi_{XY}(\mathcal{S})$, the projection of $\mathcal{S}$ on $X\times Y$, is a maximizing set for \eqref{ASSOC-2MARG-TO-3-MAR-QUAD-LAYOUT},  where $\pi_{XY}$ denotes the projection map $\pi_{XY}:X\times Y\times Z\to X\times Y $. Moreover, the potentials $\phi_1$ and $\phi_2$ solve its dual. We shall show that the maximizing sets $\mathcal{S}$ and $\pi_{XY}(\mathcal{S})$ are $(c,Q)$- and $(c_1,P)$-extreme, respectively, and then the claim follows from Theorem \ref{TH-UNI-MULT-CP-EXT}. Here 
\begin{align*}
c_1(x,y)=\langle x,y\rangle+\xi(x+y),\quad P=\{Y_k\}_{k=1}^K,\quad \text{and}\quad Q=\{Z_l\}_{l=1}^L,
\end{align*}
where
\begin{align*}
    Y_k=\overline{Y}_k\times\{\overline{y}_k\},\ k=1,\ldots,K\quad\text{and}\quad Z_l=\overline{Z}_l\times\{\overline{z}_l\},\ l=1,\ldots,L.
\end{align*}
\begin{itemize}
    \item 
\textbf{$c$-extremality of $\mathcal{S}$.}

We begin with the problem \eqref{3-MAR-QUAD-LAYOUT}. Define the following set-valued functions, 
\begin{align*}
    &F:X\times Y\to 2^Z,\quad F(x,y)=\left\{z\in Z\ :\ (x,y,z)\in\mathcal{S}\right\},\\
    &\iota : \text{Dom}(F)\to\{1,\ldots,L\},\quad \iota(x,y)=\min\{l\in\{1,\ldots,L\}\ :\ F(x,y)\cap Z_l\neq\emptyset \},\\
    &f_Q:X\times Y\to 2^Z,\quad f_Q(x,y)=\argmax\{c(x,y,z)\ : z\in F(x,y)\cap Z_{\iota(x,y)}\}.
\end{align*}
We verify the following conditions,
\begin{enumerate}
    \item $\text{Dom}(F)\cap M\times N=\text{Dom}(f_Q)\cap M\times N$.
    \item For all $(x_1,y_1),\, (x_2,y_2)\in \text{Dom}(F)\cap M\times N $, we have
    \begin{align*}
        F(x_1,y_1)\setminus\{z_1\}\cap F(x_2,y_2)\setminus\{z_2\}\cap W=\emptyset, \quad \forall z_i\in f_Q(x_i,y_i),\ i=1,2 .
    \end{align*}
\end{enumerate}
    Here, the subsets $M\subseteq X$, $N\subseteq Y$, and $W\subseteq Z$, are $\mu$-, $\nu$-, and $\gamma$-full measure subsets, respectively, and they should be determined. We note that the potentials are locally Lipschitz and due to the fact that $\mu$, $\nu_k$'s, for $k=2,\ldots,K$, and $\gamma_l$, for $l=2,\ldots, L$, are absolutely continuous with respect to the Lebesgue measure on $\mathbb{R}^n$, the following subsets are the ones that we are looking for
\begin{align*}
    &M=\left\{x\in X\ :\ \nabla\phi_1(x)\ \text{exists}\right\},\\
    &N=Y_1\cup\left(\bigcup_{k=2}^K\left\{y\in Y_k\ :\ \nabla\phi_2(y)\ \text{exists}\right\}\right),\\
    &W=Z_1\cup\left(\bigcup_{l=2}^L\left\{z\in Z_l\ :\ \nabla\phi_3(z)\ \text{exists}\right\}\right).
\end{align*}
Since the sets $X$ and $Y$ are compact and the function $c$ is continuous the first condition is automatically satisfied. For the second condition, we note that if there exist distinct pairs $(x_1,y_1)$ and $(x_2,y_2)$ in $\text{Dom}(F)\cap M\times N$, and the elements $z_1$ and $z_2$ in $f_Q(x_1,y_1)$ and $f(x_2,y_2)$ such that the intersection has the element $z^*$, then we have
\begin{align*}
    (x_1,y_1,z_1),\, (x_1,y_1,z^*),\, (x_2,y_2,z_2),\, (x_2,y_2,z^*)\in \mathcal{S}.
\end{align*}
Now, based on the layer $Z_l$ that contains $z^*$, we reach to a contradiction.
\begin{enumerate}
            \item[Case 1:] If $z^*\in Z_1$, then since the function $f_Q$ maps each element $x\in M$ to a subset of a unique layer $Z_l$ that intersects $F(x)$ with the smallest index $l=\iota(x)$, we have that $z_1,z_2\in Z_1$. Moreover, since $(x_1,y_1,z_1),\;(x_1,y_1,z^*)\in\mathcal{S}$, and through differentiability of $\phi_1$ at $x_1\in M$, it can be concluded that
\begin{align}\label{Eq28}
e_{n+1}=\alpha(z^*-z_1),\quad\text{for some}\;\alpha\in\mathbb{R}\backslash\{0\}.
\end{align}
Now, by multiplying both sides of \eqref{Eq28} by  $z^*-z_1$, we have that
\begin{align*}
    0= \langle e_{n+1}, z^*-z_1\rangle=\alpha\|z^*-z_1\|^2\neq 0,
\end{align*}
which  is a contradiction.
 \item[Case 2:] If $z^*\in\bigcup_{l=2}^L (Z_l\cap W)$, then by the definition of $W$, we have differentiability of $\phi_3$ at $z^*$. Therefore, from $(x_1,y_1,z^*),\;(x_2,y_2,z^*)\in\mathcal{S}$, we have
\begin{align}\label{Eq29}
e_{n+1}=\alpha(x_1-x_2+y_1-y_2),\quad\text{for some}\;\alpha\in\mathbb{R}\backslash\{0\}.
\end{align}
Multiplying both sides of \eqref{Eq29} by $x_1-x_2+y_1-y_2$, we obtain that
\begin{align}
    0= \langle e_{n+1},x_1-x_2+y_1-y_2\rangle=\alpha\| x_1-x_2+y_1-y_2\|^2.
\end{align}
Consequently,
\begin{align}\label{Eq30}
    x_1+y_1=x_2+y_2.
\end{align}
On the other hand, by $c$-monotonicity of $\mathcal{S}$, and the fact $(x_1,y_1,z^*),\;(x_2,y_2,z^*)\in\mathcal{S}$, it is obtained that
\begin{align}
    \langle x_1-x_2,y_1-y_2\rangle \geq 0.
\end{align}
Consequently, by multiplying both sides of \eqref{Eq29} by $x_1-x_2$, we have that
\begin{align*}
      0= \langle e_{n+1},x_1-x_2\rangle=\alpha\left(\| x_1-x_2\|^2+ \langle x_1-x_2,y_1-y_2\rangle\right).
\end{align*}
Therefore, $  x_1=x_2$, and together with \eqref{Eq30}, we get that $ y_1=y_2$. This is a contradiction as $(x_1,y_1)$ and $(x_2,y_2)$ are two distinct elements of  $\text{Dom}(F)\cap M\times N$.
\end{enumerate}
This completes the proof of $c$-extremality of $\mathcal{S}$.
\item \textbf{$c_1$-extremality of $\pi_{XY}(\mathcal{S})$.}

In order to prove the claim, we define the following set-valued functions,
\begin{align*}
    &G:X\to 2^Y,\quad G(x)=\left\{y\in Y\ :\ (x,y)\in\pi_{XY}(\mathcal{S})\right\},\\
    &\kappa : \text{Dom}(G)\to\{1,\ldots,K\},\quad \kappa(x)=\min\{k\in\{1,\ldots,K\}\ :\ G(x)\cap Y_k\neq\emptyset \},\\
    &f_P:X \to 2^Y,\quad f_P(x)=\argmax\{c_1(x,y)\ :\  y\in G(x)\cap Y_{\kappa(x)}\}.
\end{align*}
We shall show the following conditions hold,
\begin{enumerate}
    \item $\text{Dom}(G)\cap M=\text{Dom}(f_P)\cap M$.
    \item For all $x_1,\, x_2\in \text{Dom}(G)\cap M$, we have
    \begin{align*}
        G(x_1)\setminus\{y_1\}\cap G(x_2)\setminus\{y_1\}\cap N=\emptyset, \quad \forall y_i\in f_P(x_i),\ i=1,2 .
    \end{align*}
\end{enumerate}
    Here, the subsets $M\subseteq X$ and  $N\subseteq Y$ are  $\mu$-, and $\nu$-full measure subsets from the previous step. The first condition holds true via compactness of $Y$ and continuity of the function $c_1$. To verify the second condition, suppose the contrary, and let $x_1$ and $x_2$ be two distinct elements in $\text{Dom}(G)\cap M$, such that there exist $y_i\in f_P(x_i)$ for $i=1,2$, where the intersection includes the element $y^*$. Hence,
\begin{align*}
    (x_1,y_1),\, (x_1,y^*),\, (x_2,y_2),\, (x_2,y^*)\in \pi_{XY}(\mathcal{S}).
\end{align*}
Now, based on the position of $y^*$, we derive a contradiction. Before going further, we note that wherever the function $c_1$ is differentiable with respect to $x$ or $y$, we have that,
\begin{align*}
    \nabla_xc(x,y)=y+\nabla\xi(x+y),\quad\text{and}\quad \nabla_yc(x,y)=x+\nabla\xi(x+y).
\end{align*}
    \begin{enumerate}
            \item[Case 1:] $y^*\in Y_1$. By the definition of $f_P$ that relies on the order of the layers $Y_k$, we have that $y_1,y_2\in Y_1$. Additionally, by differentiability of $\phi_1$ at $x_1$ and $(x_1,y_1),\;(x_1,y^*)\in\pi_{XY}(\mathcal{S})$, we conclude that
\begin{align}\label{Eq31}
e_{n+1}=\alpha(y^*-y_1+\nabla\xi(x_1+y^*)-\nabla\xi(x_1+y_1) ),\quad\text{for some}\;\alpha\in\mathbb{R}\backslash\{0\}.
\end{align}
Now, by multiplying both sides of \eqref{Eq31} by  $y^*-y_1$, and through convexity of $\xi$, we have that
\begin{align*}
    0= \langle e_{n+1}, y^*-y_1\rangle=\alpha\left(\|y^*-y_1\|^2+\langle \nabla\xi(x_1+y^*)-\nabla\xi(x_1+y_1), y^*-y_1\rangle\right)\neq 0,
\end{align*}
which  is a contradiction.
 \item[Case 2:] $y^*\in\bigcup_{k=2}^K (Y_k\cap N)$. In this case we can differentiate the potential $\phi_2$ at the point $y^*$ with respect to any coordinate chart representing the $Y_k$ that contains $y^*$. Therefore, by the fact that $(x_1,y^*),\;(x_2,y^*)\in \pi_{XY}(\mathcal{S})$, we have
\begin{align}\label{Eq32}
e_{n+1}=\alpha(x_1-x_2+\nabla\xi(x_1+y^*)-\nabla\xi(x_2+y^*)),\quad\text{for some}\;\alpha\in\mathbb{R}\backslash\{0\}.
\end{align}
Multiplying both sides of \eqref{Eq32} by $x_1-x_2$, and relying on monotonicity of the derivative of the convex function $\xi$, we obtain that
\begin{align}
    0= \langle e_{n+1},x_1-x_2\rangle=\alpha\left(\| x_1-x_2\|^2+\langle \nabla\xi(x_1+y^*)-\nabla\xi(x_2+y^*), x_1-x_2\rangle \right)\neq 0,
\end{align}
which  is a contradiction. Thus, the set $\pi_{XY}(\mathcal{S})$ is $c_1$-extreme.
\end{enumerate}
    \end{itemize}
Now, by applying Theorem \ref{TH-UNI-MULT-CP-EXT} we obtain that the solution to \eqref{3-MAR-QUAD-LAYOUT} is unique. In the remainder of the proof, we shall demonstrate that the unique solution $\lambda$ is concentrated on the union of the graph of functions $T_{k}:X\to Y\times Z$, where $k\in \{1,\ldots, K\times L\}$. To achieve this, we will establish that for any $(x_0,y_0,z_0)\in M\times N\times W$, the following set has at most $K\times L$ elements, and then we will apply Theorem \ref{M-TWIST-TH}
\[ L(x_0,y_0,z_0)=\left\{(y,z)\in (Y\times Z)\cap (N\times W) \ :\ \nabla_xc(x_0,y,z)=\nabla_xc(x_0,y_0,z_0) \right\}. \]

Let it be otherwise and suppose there are $K\times L +1$ points $\{(y_j,z_j)\}_{j=1}^{K\times L+1}$ in the aforementioned set. Then, there exists at least one $(k_0,l_0)\in\{1,\ldots,K\}\times \{1,\ldots,L\}$ that contains two elements, $(y_{j_1},z_{j_1})$ and $(y_{j_2},z_{j_2})$, that is, $y_{j_k}\in Y_{k_0}$ and $z_{j_l}\in Z_{l_0}$, for $k,l=1,2$. Consequently,
\[ (x_0,y_{j_1},z_{j_1}),\, (x_0,y_{j_2},z_{j_2})\in\mathcal{S}. \]
Through the differentiability of $\phi_1$ at $x_0\in M$, we derive
  \begin{align}\label{Eq33}
      e_{n+1}=\alpha(y_{j_1}-y_{j_2}+z_{j_1}-z_{j_2}),\quad\text{for some}\;\alpha\in\mathbb{R}\backslash\{0\}.
  \end{align}
Additionally, employing the $c$-monotonicity of $\mathcal{S}$, we deduce
  \begin{align*}
      \langle z_{j_1}-z_{j_2},y_{j_1}-y_{j_2}\rangle\geq 0.
  \end{align*}
  Consequently, by multiplying both sides of \eqref{Eq33} by $y_{j_1}-y_{j_2}$, we find
  \begin{align*}
          0=\langle e_{n+1},y_{j_1}-y_{j_2}\rangle =\alpha\left(\|y_{j_1}-y_{j_2}\|^2+\langle z_{j_1}-z_{j_2}, y_{j_1}-y_{j_2}\rangle\right)\neq 0,
  \end{align*}
  which leads to $y_{j_1}=y_{j_2}$. Similarly, applying the same process, we conclude that $z_{j_1}=z_{j_2}$. Here, we used the facts that $y_{j_1}-y_{j_2}$ and $z_{j_1}-z_{j_2}$ are orthogonal to $e_{n+1}$. This contradicts the assumption of $(y_{j_1},z_{j_1})$ and $(y_{j_2},z_{j_2})$ being distinct elements. So, the cost function $c$ satisfies the $m$-twist condition for $m=K\times L$ and the unique optimal plan $\lambda$ is concentrated on the graph of $K\times L$ functions, by Theorem \ref{M-TWIST-TH}.
     \end{proof}

\end{theorem}

\begin{remark}
We would like to highlight that when $K=L=1$, the result of the preceding theorem aligns with those established by Gangbo and \'Swi\c ech \cite{GANGBOSWIECH} for multi-marginal optimal mass transportation. Specifically, they investigated the $N$-marginal Monge-Kantorovich problem with a quadratic cost function, demonstrating uniqueness under the assumption of absolute continuity of the $N$ marginals. An alternative proof is also presented in \cite{AH-MOMENI} through a reduction argument.
\end{remark}

\section{Cost functions with the sub-twist property}\label{SECTIONSUBTWIST}
The problem \eqref{TAGMKP} with a continuous cost function $c:X\times Y\to\mathbb{R}$ defined on the compact subsets $X$ and $Y$ in $\mathbb{R}^{n+1}$ that satisfies twist, sub-twist, $m$-twist, or generalized twist conditions  has been studied extensively \cite{AHMADKIMMCCANN,AMBROSIOKIRCHHEIMPRATELLI,BERNARDBORIS,CAFFARELLI,CARLIER,CHAMPIONDEPASCALE2,CHAMPIONDEPASCALE,CHAMPIONDEPASCALE3,FATHIFIGALLI,FIGALLI22,GANGBOMCCANN,LEVIN,MCCANN222,MOMENI-CHAR}. It has been shown that, under absolutely continuity of the first marginal with respect to the Lebesgue measure on $\mathbb{R}^{n+1}$,  \eqref{TAGMKP} admits a unique solution which is concentrated on, the graph of a function $T:X\to Y$, the union of the graphs of two functions $T:X\to Y$ and $S:Y\to X$, the union of the graph of $m$ functions $T_k:X\to Y$, and on the union of countable number of graphs of functions $T_k:X\to Y$, respectively.

In this section, we consider a two-marginal optimal transport problem corresponding to a cost function \( c: X \times Y \to \mathbb{R} \) which satisfies the celebrated sub-twist condition introduced in \cite{CHIAPPORIMCCANNNESHEIM}. Our goal is to study the uniqueness of the optimal plan for the two-marginal Monge-Kantorovich problem on \( X \times Y \subseteq \mathbb{R}^{n+1} \times \mathbb{R}^{n+1} \) in the case where the probability measure on \( X \) assigns positive measures to some \( n \)-dimensional sub-manifold \( X_0 \) of \( X \). More precisely, we consider compact subsets \( X, Y \subseteq \mathbb{R}^{n+1} \) with smooth boundaries, and equip \( X \) and \( Y \) with the probability measures \( \mu \) and \( \nu \), where \( \mu \) is defined as follows
\begin{align*}\tag{FM}\label{TAGFM}
    \int_X f(x)\,d\mu(x)=\int_X f(x)\alpha(x)\,d\mathcal{L}^{n+1}(x)+\int_{X_0} f(x_0)\beta(x_0)\, d\mathcal{L}^{n}(x_0),\quad \forall f\in C_b(X).
\end{align*}
Here, \( d\mathcal{L}^{n+1} \) and \( d\mathcal{L}^{n} \) stand for the Lebesgue measures on \( \mathbb{R}^{n+1} \) and \( \mathbb{R}^{n} \), respectively. Moreover, \( X_0 \) is any \( n \)-dimensional sub-manifold of \( X \) such that either \( X_0 \) has no boundary, or \( \partial X_0 \), the boundary of \( X_0 \), is a smooth \( (n-1) \)-dimensional manifold. Additionally, \( \alpha \) and \( \beta \) are two non-negative measurable functions residing in \( L^\infty(X) \) and \( L^\infty(X_0) \), respectively. For instance, the subset \( X_0 \) can be chosen to be the boundary of the set \( X \). The first marginal \( \mu \) can be expressed as
\begin{align}\label{DECOPM-MU}
    \mu=s\mu_1+(1-s)\mu_2,
\end{align}
where
\begin{align*}
    d\mu_1=\frac{1}{s}\alpha d\mathcal{L}^n,\quad     d\mu_2=\frac{1}{1-s}\beta d\mathcal{L}^{n-1},\quad\text{and}\quad s=\mu(X\setminus X_0).
\end{align*}
The main feature distinguishing the current problem from previous works is the second term in the definition of the measure \( \mu \) in  \eqref{TAGFM}, namely, the measure \( \mu_2 \) assigning measure to the \( n \)-dimensional manifold. To state our first result, we need to recall the following definition.
\begin{definition}
    We say that a function $c:X\times Y\to\mathbb{R}$ satisfies the sub-twist condition when for each pair of distinct elements $y_1,y_2\in Y$ the function 
    \begin{align*}
x\mapsto c(x,y_1)-c(x,y_2),
    \end{align*}
    has no critical point except one global minimum and one global maximum.
\end{definition}

Before going further, let us elaborate more about critical points of a function $H:X\to \mathbb{R}$ where $X$ can be an $n$-dimensional manifold in $\mathbb{R}^{n+1}$. Here is the definition of a critical point of a differentiable function on a differentiable manifold.
\begin{definition}\label{DEF-CRITICAL}
    Let $X\subseteq \mathbb{R}^{n+1}$ be an $n$-dimensional differentiable manifold, and let $H:\mathbb{R}^{n+1}\to \mathbb{R}$ be a differentiable map in the usual sense. We say that $x_0\in X$ is a critical point of $H|_{X}$, if there exists a scalar $\alpha\in\mathbb{R}\setminus\{0\}$ such that
    \begin{align*}
\nabla H(x_0)=\alpha\mathrm{n}(x_0),
    \end{align*}
    where $\mathrm{n}(x_0)$ denotes the normal vector to $X$ at the point $x_0$.
\end{definition}

Now we state our  result in the following theorem.
\begin{theorem}\label{TH-2MARG-SUBTWIST}
Let $(X,\mathcal{B}_X,\mu)$ and $(Y,\mathcal{B}_Y,\nu)$ be two Borel probability spaces where $X,Y\subseteq\mathbb{R}^{n+1}$ are compact subsets, and $X_0$ be any compact sub-manifold of $X$ of $n$-dimension that either $X_0$ has no boundary, or, $\partial X_0$, the boundary of $X_0$ is a smooth $(n-1)$-dimensional manifold. Assume that $\mu$ is a Borel probability measure of the form \eqref{TAGFM}. Let $c:X\times Y\to \mathbb{R}$ be a differentiable map that satisfies the sub-twist condition. Then the optimization problem
\begin{align}\label{NEW-TWO-MARG}
    \inf\left\{\int_{X\times Y}c(x,y)d\,\lambda\ :\ \lambda\in\Pi(\mu,\nu)\right\},
\end{align}
admits a unique solution.
\begin{proof}
Clearly, solutions to the problem \eqref{TH-2MARG-SUBTWIST} and its dual exist. Assume that $(\phi,\psi)\in L_1(X,\mu)\times L_1(Y,\nu)$ are potentials. To prove uniqueness, we shall show that the minimizing set $\mathcal{S}$ of \eqref{NEW-TWO-MARG} is $c$-extreme. Recall the following set and set-valued maps,
\begin{align*}
&\mathcal{S}= \{(x,y)\in X\times Y\ :\ c(x,y)=\phi(x)+\psi(y) \},\\
   & F:X\to 2^Y,\quad F(x)=\{ y\in Y\ :\ (x,y)\in\mathcal{S}\},\\
   &f:X\to 2^Y,\quad f(x)=\argmax\{c(x,y)\ :\ y\in F(x)\}.
\end{align*}
We shall show that the two conditions in Definition \ref{DEF-C-EXT} are satisfied for two $\mu$- and $\nu$-full measure subsets $M$ and $N$ of $X$ and $Y$, respectively. Due to the compactness of $X$, condition $(i)$ in Definition \ref{DEF-C-EXT} is already satisfied for any $\mu$-full measure subset $M\subseteq X$. To show the second condition, we need to determine the $\mu$- and $\nu$-full measure subsets $M\subseteq X$ and $N\subseteq Y$. According to Lemma \ref{LEMMAEXISTENCE}, the potentials $\phi$ and $\psi$ are locally Lipschitz. Define
\begin{align*}
    M_1:=\left\{x\in X\setminus X_0\ :\ \nabla\phi(x)\;\text{exists} \right\}\quad\text{and}\quad     M_2:=\left\{x\in X_0\ :\ \nabla\phi(x)\;\text{exists} \right\}
\end{align*}
In the set $M_2$, the derivative of $\phi$ is considered as the derivative on a manifold. Note that since the potential $\phi$ is locally Lipschitz and the measures $\mu_1$ and $\mu_2$ are absolutely continuous with respect to $\mathcal{L}^{n+1}$ and $\mathcal{L}^{n}$ on each coordinate chart of $X_0$, respectively, it follows that the sets $M_1$ and $M_2$ are $\mu_1$- and $\mu_2$-full measure sets, where $\mu_i$'s are measures described in \eqref{DECOPM-MU}. Consequently, the subset
\begin{align*}
    M_0=M_1\cup M_2,
\end{align*}
is a $\mu$-full measure subset of $X$. The subset $N$ can be chosen arbitrarily. Without loss of generality, we assume that there exists $\overline{x}\in X$ such that $c(\overline{x},y)=0$ for all $y\in Y$, simply by replacing $c$ with the function $\overline{c}(x,y):=c(x,y)-c(\overline{x},y)$. Let $M=M_0\setminus\{\overline{x}\}$, and obviously, $\mu(M)=1$. We shall show that for distinct elements $x_1,x_2\in M$, we have
\begin{align}\label{RR}
   \big( F(x_1)\setminus\{y_1\}\big)\cap  F(x_2)\cap N=\emptyset,\quad \forall y_1\in f(x_1).
\end{align}
To this end, we show that for each $x_1\in M$, the set $f(x_1)$ is a singleton. Suppose otherwise, and let there be two distinct elements $y_1$ and $y_1'$ in $f(x_1)$. Then, $x_1$ is a critical point of the function
\begin{align*}
    x\mapsto c(x,y_1)-c(x,y_1').
\end{align*}
Since, the function saves at most one global minimum and one global maximum, it implies that
\begin{align*}
    c(x,y_1)-c(x,y_1')\neq 0,\quad \forall x_1\in M \ \text{with}\ x\neq x_1. 
\end{align*}
But this is a contradiction as $c(\overline{x},y_1)-c(\overline{x},y_1')=0$. Thus, $f(x_1)$ is singleton. Let us denote its element by $y_1$, i.e., $f(x_1)=\{y_1\}$. Now, we prove \eqref{RR}. If there exists an element $y^*$ in this intersection, then $x_1$ is a critical point of the function
\begin{align}\label{SUBMAP}
x\mapsto c(x,y^*)-c(x,y_1).
\end{align}
Furthermore, since $f(x_1)=\{y_1\}$, we have
\begin{align*}
    c(x_1,y_1)>c(x_1,y^*),
\end{align*}
implying that $x_1$ is a global minimum of the map \eqref{SUBMAP}. Hence, we obtain
\begin{align*}
    c(x_2,y^*)-c(x_2,y_1)> c(x_1,y^*)-c(x_1,y_1),
\end{align*}
which we can rearrange as
\begin{align}\label{CONTA}
    c(x_1,y_1)+c(x_2,y^*)> c(x_2,y_1)+c(x_1,y^*).
\end{align}
On the other hand, we note that 
\begin{align*}
    (x_1,y_1),\ (x_2,y^*)\in \mathcal{S}.
\end{align*}
Thus, by the $c$-cyclical monotonicity of $\mathcal{S}$, we have
\begin{align*}
    c(x_1,y_1)+c(x_2,y^*)\leq c(x_2,y_1)+c(x_1,y^*).
\end{align*}
This last inequality contradicts \eqref{CONTA}. Therefore, the intersection \eqref{RR} is empty, and $\mathcal{S}$ is a $c$-extreme minimizing set. Consequently, the solution to the optimization problem \eqref{NEW-TWO-MARG} is unique.

\end{proof}
\end{theorem}

\section{Applications}\label{SECTIONEXAMPLES}
In this section, we investigate the uniqueness of optimal plans for two-marginal Monge-Kantorovich problems, specifically focusing on cases where the quadratic cost function is employed alongside a first marginal structured as defined in \eqref{TAGFM} where the $n$-dimensional subset $X_0$ of $X$ is considered to be the boundary of $X$, that is $X_0=\partial X$. Broadly viewed, this investigation can be construed as direct outcomes of findings presented in the preceding sections. By delving into this specialized application, we aim to demonstrate the the particular structural constraints imposed by the form of the first marginal, helping us understand the special behaviors of the unique optimal plan that appear in this situation.

\noindent Consider $X,Y\subseteq\mathbb{R}^{n+1}$ where $X$ is a compact subset with the smooth boundary and $Y$ is any compact subset. Equip $X$ and $Y$ with the Borel probability measures $\mu$ and $\nu$ where $\mu $ is defined as follows
\begin{align}\label{Q1}
    \int_X f(x)\,d\mu(x)=\int_X f(x)\alpha(x)\,d\mathcal{L}^{n+1}(x)+\int_{\partial X} f(\overline{x})\beta(\overline{x})\, dS(\overline{x}),\quad \forall f\in C_b(X),
\end{align}   
and $\nu$ is any Borel probability measure. Here, $dS$ is the surface area on $\partial X$. Our goal is to fully determine the optimal plans of the following optimal transport problem
\begin{align}\label{Q2}
    \inf\left\{\int_{X\times Y}|x-y|^2d\lambda\ :\  \lambda\in \Pi(\mu,\nu)\right\}.
\end{align}

\begin{theorem}\label{FIRSTTH}
    Let $X,Y\subseteq\mathbb{R}^{n+1}$ be such that $X$ and $Y$ are compact domains, $X$ has the smooth boundary equipped with the Borel probability measures $\mu$ as \eqref{Q1} and $Y$ is associated with any Borel probability measure $\nu$. Then any optimal plan $\lambda\in \Pi(\mu,\nu)$ of \eqref{Q2} transports $\mu$-almost every points in $\text{int}(X)$, interior of $X$,  to a single point while it transports each  point in $\partial X$ to a set that is a subset of a  straight line.    

\begin{proof}
First, analogous to\eqref{DECOPM-MU}, we observe that
\begin{align*}
    \mu = s\mu_1 + (1-s)\mu_2, 
\end{align*}
where the two probability measures $\mu_1$ and $\mu_2$ are normalized restrictions of $\mu$ to $\text{int}(X)$ and $\partial X$, respectively, and $s=\mu(\text{int}(X))$.
Since the domains are compact and the cost function is continuous, then by Lemma \ref{LEMMAEXISTENCE}, the problem \eqref{Q2} and its dual admit solutions. For locally Lipschitz potentials $\phi$ and $\psi$ on $X$ and $Y$, we have
\begin{align*}
    c(x,y) \geq \phi(x) + \psi(y),\quad \forall (x,y)\in X\times Y,
\end{align*}
and consider the following minimizing set
\begin{align*}
   \mathcal{S} = \{(x,y)\in X\times Y\ :\ c(x,y) = \phi(x) + \psi(y)\},
\end{align*}
which has $\lambda$-full measure for every optimal plan $\lambda$ of \eqref{Q2}.  Let
\begin{align*}
    M_1 = \{x\in \text{int}(X) \ :\ \nabla\phi(x)\; \text{exists}\}, \quad\text{and}\quad M_2 = \{\overline{x}\in \partial X  \ :\ \nabla\phi(\overline{x})\; \text{exists}\}.
\end{align*}
Since $\phi$ is locally Lipschitz and $\mu_1$ is absolutely continuous with respect to $\mathcal{L}^{n+1}$, the set $M_1$ has $\mu_1$-full measure. Similarly, the subset $M_2$ has $\mu_2$-full measure. Therefore, the set $M=M_1\cup M_2$ has $\mu$-full measure. Now, by applying Brenier's result on $\text{int}(X)$, we find out that there exists a unique optimal transport map $T:\text{int}(X)\to Y$ such that $\lambda|_{\text{int}(X)\times Y} = (\text{id}\times T)\#\mu$, meaning that elements in $\text{int}(X)$ are uniquely transported to elements in $Y$. To handle elements in $\partial X$, let $(\overline{x},y_1)$ and $(\overline{x},y_2)$ be two elements in $\mathcal{S}$ such that $\overline{x}\in\partial X$. Then via Lemma \ref{LEM-CRITICAL}, we obtain that
\begin{align*}
    y_2 = \alpha\mathrm{n}(\overline{x}) + y_1,\quad \text{for some}\ \alpha\in\mathbb{R}\setminus\{0\}.
\end{align*}
This reveals that any element in $F(\overline{x}) = \{y\in Y\ :\ (\overline{x},y)\in\mathcal{S}\}$ lies on the line parallel to $\mathrm{n}(\overline{x})$ and passing through $y_1$.

\end{proof}

\end{theorem}


\begin{corollary}\label{CORO66}
    If in Theorem \ref{FIRSTTH} the subset $X$ is a strictly convex bounded subset  of $\mathbb{R}^{n+1}$ with the smooth boundary, then the optimal transport problem \eqref{Q2} admits a unique solution.

    \begin{proof}
        To show the uniqueness, we shall prove that an optimal plan of \eqref{Q2} is $c$-extreme. To do so, consider the two set-valued functions,
        \begin{align*}
            &G:Y\to 2^X,\quad G(y)=\{x\in X\ :\ (x,y)\in \mathcal{S}\},\\
            &g:Y\to 2^X.\quad g(y)=\argmax\{c(x,y)\ : x\in G(y)\}.
        \end{align*}  
        We must show that the following criteria are satisfied
\begin{enumerate}
    \item  $\text{Dom}(G)\cap N=\text{Dom}(g)\cap N$,
    \item $G(y_1)\backslash \{x_1\}\cap G(y_2)\backslash \{x_2\}\cap M=\emptyset$, for $\nu$-a.e. $y_i\in Y$ with $y_1\neq y_2$ and all $x_i\in g(y_i)$.
\end{enumerate}
Here, the subset $M$ is given in the proof of Theorem \ref{FIRSTTH}, and $N$ can be chosen to be any $\nu$-full measure subset of $Y$. Remark that the first condition   is satisfied by the compactness of $X$ and   continuity of $c$. To show the second one, we assume otherwise. For any element $x^*$ in the intersection, we have the following cases based on the position of $x^*$.
        \begin{enumerate}
            \item[Case 1:] If $x^*\in \text{int}(X)\cap M$, then by differentiability of $\phi$ at $x^*$, in the usual sense, and the fact that $(x^*,y_1),(x^*,y_2)\in\mathcal{S}$, we obtain that
\begin{align*}
    2x^*-2y_1=2x^*-2y_2,
\end{align*}
which implies that $y_1=y_2$, that is a contradiction.
\item[Case 2:] If $x^*\in \partial X\cap M$, then by smoothness of the boundary, and $(x^*,y_1),(x^*,y_2)\in\mathcal{S}$, we obtain that
\begin{align}\label{Q3}
    \mathrm{n}(x^*)=\alpha(y_2-y_1),\quad\text{for some}\;\alpha\in\mathbb{R}\setminus\{0\}.
\end{align}
On the other hand, by strict convexity of the domain $X$, we have
\begin{align}\label{Q4}
    \langle \mathrm{n}(x^*),x^*-x_1\rangle >0,\qquad \langle \mathrm{n}(x^*),x^*-x_2\rangle >0.
\end{align}
Additionally, since $\{(x_1,y_1),(x^*,y_2)\},\; \{(x^*,y_1),(x_2,y_2)\}\subseteq\mathcal{S}$, then via $c$-cyclic monotonicity of $\mathcal{S}$, we obtain that
\begin{align}\label{Q5}
    \langle x^*-x_1,y_2-y_1 \rangle\geq 0,\qquad    \langle x^*-x_2,y_1-y_2 \rangle\geq 0,
\end{align}
By substituting \eqref{Q3} into \eqref{Q4} and comparing it to \eqref{Q5}, we conclude that $\alpha$ must be both positive and negative, which leads to a contradiction.
        \end{enumerate}
Thus, the minimizing set $\mathcal{S}$ is $c$-extreme. Consequently, by Theorem \ref{TH-C-EXT}, the optimal plan of \eqref{Q2} is unique.        
    \end{proof}
\end{corollary}

\begin{remark}
    In the proof of Corollary \ref{CORO66}, we utilized set-valued functions with domains as subsets of $Y$ and with the range into $2^X$, which indicates the flexibility of the notion of $c$-extremal sets.
\end{remark}

\section{Appendix}\label{SECTIONAPPENDIX}
In this section, we provide details regarding some technical justifications of the methods implemented in this paper. Primarily, one might be concerned about the differentiability of the potentials and some generalizations concerning minimizing sets for the $N$-marginal case, which holds when $N=2$.

\subsection{Differentiability of the potentials}
The forthcoming lemma addresses an essential aspect of the dual problem's solutions: their differentiability. Understanding the differentiability properties of these potentials is crucial for analyzing the behavior of optimal transport maps and gaining insights into the underlying structures of the problem. This lemma sheds light on the conditions under which the potentials arising from the dual problem exhibit differentiability, providing a deeper understanding of their mathematical properties and implications for optimization strategies. Here, the domain $X\subseteq\mathbb{R}^{n+1}$ is a differentiable manifold of dimension $n$. This implies that for each $x\in X$, any chart $(f_x,U_x)$ representing it on $B_r(0)\subseteq\mathbb{R}^n$ is differentiable, that is, the map $f_x:B_r(0)\to U_x$, with $f(0)=x$ is differentiable: There exist a linear operator $Df(0)$ such that
\begin{align*}
    \lim_{h\to 0}\frac{|f(h)-f(0)-\langle Df(0),h\rangle|}{|h|}=0.
\end{align*}
Here $U_x\subseteq X$ is a relatively open subset containing $x$, i.e., $U_x=X\cap U$, for some open set $U\subseteq \mathbb{R}^{n+1}$.
\begin{lemma}\label{LEM-DIFF}
    Let $X,Y\subseteq\mathbb{R}^{n+1}$ be two compact $n$-dimensional differentiable manifolds which either have no boundary or their boundaries are differentiable $(n-1)$-dimensional manifolds equipped with  the Borel probability measures $\mu$ and $\nu$, respectively. Let $c:X\times Y\to\mathbb{R}$ be a differentiable function. If $\mu$ is absolutely continuous with respect to $\mathcal{L}^n$ on each coordinate chart of $X$, then there exists a $c$-concave function $\phi:X\to\mathbb{R}$ (given by Lemma \ref{LEMMAEXISTENCE}) which is $\mu$-almost everywhere differentiable in a manifold sense, that is, there exists a $\mu$-full measure set $M$ such that for every $x\in M$, and for every chart $(f_x,U_x)$ the function $\phi \circ f_x$ is differentiable at $t=0$ in the usual sense. Moreover,
    \begin{align}\label{CHAINRULE}
        \nabla(\phi \circ f_x)(0)=\langle \nabla_x c(x,y), Df_x(0)\rangle,\quad\text{for some $y\in F(x)$} .
    \end{align}
    where the set-valued function $F$ is introduced in \eqref{F-function}.
    \begin{proof}
       By Lemma \ref{LEMMAEXISTENCE} existence of such a $c$-concave potential $\phi$ which is locally Lipschitz (with respect to the metric on $\mathbb{R}^{n+1}$) is guaranteed. Therefore, it is almost everywhere differentiable on $\mathbb{R}^{n+1}$. Thus, due to absolute continuity of $\mu$ with respect to $\mathcal{L}^{n}$, the function $\phi|_X$ is differentiable over $X$. This means that there exits a subset $M\subseteq X$ with $\mu(M)=1$, and is such that for any $x\in M$, and any differentiable chart $(f_x,U_x)$, the function $\phi|_X\circ f_x $ is differentiable at $t=0\in\mathbb{R}^n$, in the usual sense. On the other hand, the subset $M$ can be chosen to be a subset of $\pi_X(\mathcal{S})$, where $\mathcal{S}$ is a minimizing set for the problem \eqref{TAGMKP} and $\pi_X$ denotes projection map from $X\times Y$ onto $X$. Therefore, due to existence of $\nabla_xc(x,y)$, for some $y\in F(x)$, we have \eqref{CHAINRULE}.
    \end{proof}
\end{lemma}

By utilizing Lemma \ref{LEM-DIFF}, we immediately conclude the following observation which we state it as a lemma.

\begin{lemma}\label{LEM-CRITICAL}
Under the assumptions of Lemma \ref{LEM-DIFF}, there exists a $\mu$-full measure subset $M\subseteq X$ such that for every $x\in M$, if $y_i\in F(x)$, for $i=1,2$, then we have
    \begin{align}\label{NORMALPARALLEL}
        \nabla_x c(x,y_1)-\nabla_x c(x,y_2)=\alpha\mathrm{n}(x),\quad\text{for some}\ \alpha\in\mathbb{R}\setminus\{0\},
    \end{align}
    where $\mathrm{n}(x)$ denotes the unit outward normal to the surface $X$ at the point $x$.
    \begin{proof}
        The proof passes through the fact that for such points $(x,y_i)\in\mathcal{S}$, for $i=1,2$, from \ref{CHAINRULE}, we have that
        \begin{align*}
            \langle \nabla_x c(x,y_1)-\nabla_x c(x,y_2), Df_x(0)\rangle=0, 
        \end{align*}
        which implies \eqref{NORMALPARALLEL}.
    \end{proof}
\end{lemma}

\begin{remark}
    The reason for our approach in Definition \ref{DEF-CRITICAL} in defining critical points of  a function $H:X\to\mathbb{R}$ is in fact Lemma \ref{LEM-DIFF}. Indeed, via \eqref{CHAINRULE}, we have
    \begin{align*}
            \nabla(H|_X \circ f_x)(0)=\langle \nabla H(x), Df_x(0)\rangle . 
    \end{align*}
    Having a critical point at $x=x_0$ means the right hand side equals zero at $x=x_0$, which implies $\nabla H(x_0)=\alpha\mathrm{n}(x_0)$, for some $\alpha\in\mathbb{R}\backslash\{0\}$.
\end{remark}

\subsection{The minimizing set in $N$-marginal case}
In dealing with multi-marginal problems, one of the approaches can be done through a reduction argument. In fact, by using the solutions to the dual problem, one can associate a multi-marginal Monge-Kantorovich problem with lower-marginal problem (see \cite{AH-MOMENI}). To clarify more, for a family of Borel probability spaces $\{(X_i,\mathcal{B}_{X_i},\mu_i)\}_{i=1}^N$ with $N\geq 3$, consider the following $N$-marginal Monge-Kantorovich problem and its dual,

\begin{align*}\tag{MMKP}\label{TAGMMKP}
    &\inf\left\{\int c(x_1,\ldots,x_N)\, d\lambda\ :\ \lambda\in\Pi(\mu_1,\ldots,\mu_N) \right\},
\end{align*}
and
\begin{align*}\tag{DMMKP}\label{TAGDMMKP}
        &\sup\left\{\sum_{i=1}^N\int_{X_i} \phi_i(x)\,d\mu_i(x)\ :\ c(x_1,\ldots,x_N)\geq \sum_{i=1}^N\phi(x_i),\ \phi_i\in L_1(X_i,\mu_i) \right\}.
\end{align*}
The following theorem can be found in \cite{KELLERER} that we state it in our setting. Indeed, the following theorem partially holds in a more general setting.
\begin{theorem}\label{TH-DUALATTAIN}
    Assume that $\{(X_i,\mathcal{B}_{X_i},\mu_i)\}_{i=1}^N$ is a collection of Borel probability spaces on $\mathbb{R}^{n+1}$ such that either
    \begin{enumerate}
        \item     the domains $X_i$'s are compact, or,
        \item the measures $\mu_i$'s are compactly supported.
    \end{enumerate}
Moreover, let $c:X_1\times\cdots\times X_N\to\mathbb{R}$ be a continuous cost function. Then there is no gap between the infimum in \eqref{TAGMMKP} and the supremum in \eqref{TAGDMMKP} and these extremum values are attained at some measure $\lambda\in\Pi(\mu_1,\ldots,\mu_N)$ and continuous potentials $(\phi_1,\ldots,\phi_N)\in L_1(X_1,\mu_1)\times\cdots\times L_1(X_N,\mu_N)$. Additionally, there exists a $c$-cyclically monotone compact set $\mathcal{S}$ that contains the support of all optimal plans,
\begin{align*}
    \mathcal{S}=\left\{(x_1,\ldots,x_N)\in\prod_{i=1}^NX_i\ :\ c(x_1,\ldots,x_N)=\sum_{i=1}^N\phi_i(x_i)\right\}.
\end{align*}
\end{theorem}
The set $\mathcal{S}$ can be considered as $c$-superdifferential of the $c$-concave function,
\begin{align}\label{PHIONE}
    \phi_1(x_1)=\inf\left\{c(x_1,x_2,\ldots,x_N)-\sum_{i=2}^N\phi_i(x_i)\ :\ (x_2,\ldots,x_N)\in\prod_{i=2}^N X_i\right\}.
\end{align}
By using the potentials $\phi_i$'s one can define the following cost functions,
\begin{align}\label{C-J-COSTS}
    c_j(x_1,\ldots,x_j)=\inf\left\{ c(x_1,\ldots,x_N)-\sum_{i=j+1}^N\phi(x_i)\ :\ (x_{j+1},\ldots,x_N)\in\prod_{i=j+1}^N X_i\right\},\quad j=2,\ldots,N-1.
\end{align}
Moreover, we can associate to $N$-marginal problem \eqref{TAGMMKP} the following $j$-marginal problems,
\begin{align*}\tag{RMKP}\label{TAGRMKP}
    &\inf\left\{\int c_j(x_1,\ldots,x_j)\, d\tau\ :\ \tau\in\Pi(\mu_1,\ldots,\mu_j) \right\},\quad j=2,\ldots,N-1.
\end{align*}

Regarding the relation between optimal plans of \eqref{TAGMMKP} and \eqref{TAGRMKP} we have the following proposition. Let $\pi_j:\prod_{i=1}^N X_i\to\prod_{i=1}^j X_i$ be the project map $\pi_j(x_1,\ldots,x_N)=(x_1,\ldots,x_j)$ and for a measure $\lambda\in\Pi(\mu_1,\ldots,\mu_N)$, we denote by $\lambda_j$ the restriction of $\lambda$ on $\prod_{i=1}^j X_i$ which belongs to $\Pi(\mu_1,\ldots,\mu_j)$, that is, $\lambda_j=\pi_j\#\lambda$.
\begin{proposition}\label{PROP-OPTIMAL}
    If $\lambda\in\Pi(\mu_1,\ldots,\mu_N)$ is an optimal plan for \eqref{TAGMMKP} then $\lambda_j=\pi_j\#\lambda$ is an optimal plan for \eqref{TAGRMKP}, for $j=2,\ldots,N-1$.
    \begin{proof}
        For a proof and a more general statement see \cite[Proposition 3.1 ]{AH-MOMENI}.
    \end{proof}
  
\end{proposition}
The following lemma is one the key facts for our results in this paper.
\begin{lemma}\label{LEM-TRI-EXT}
    Let $\lambda\in\Pi(\mu_1,\ldots,\mu_N)$ and let $\lambda_j$ be the restriction of $\lambda$ to $\prod_{i=1}^jX_i$. If $\lambda$ is an extreme point of $\Pi(\lambda_{N-1},\mu_N)$ and $\lambda_j$'s are extreme points of $\Pi(\lambda_{j-1},\mu_j)$, then the measure $\lambda$ is an extreme point of $\Pi(\mu_1,\ldots,\mu_N)$.
    \begin{proof}
           If $\lambda$ is not an extreme point of $\Pi(\mu_1,\ldots,\mu_N)$, then there exist two distinct measures $\eta$ and $\zeta$ in $\Pi(\mu_1,\ldots,\mu_N)$ such that
\begin{align}\label{QQQ}
    \lambda=\frac{1}{2}(\eta+\zeta),\qquad \lambda\neq\eta ,\zeta .
\end{align}
Let $\eta_j=\pi_{j}\#\eta$ and $\zeta_j=\pi_{j}\#\zeta$ be restrictions of $\eta$ and $\zeta$ on $\prod_{i=1}^jX_i$. Then from \eqref{QQQ} we have
\begin{align*}
    \lambda_j=\frac{1}{2}(\eta_j+\zeta_j),\quad j=2,\ldots, N-1.
\end{align*}
Specially, for $j=2$, we have $ {\displaystyle \lambda_2=\frac{1}{2}(\eta_2+\zeta_2)}$. But by the assumption, $\lambda_{2}$ is an extreme point of $\Pi(\mu_1,\mu_2)$. Therefore,
\begin{align*}
    \eta_{2}=\zeta_{2}=\lambda_{2}.
\end{align*}
This implies that 
\begin{align*}
    \eta_{3},\;\zeta_{3}\in\Pi(\lambda_{2},\mu_3),
\end{align*}
and together with the facts that  $ {\displaystyle \lambda_3=\frac{1}{2}(\eta_3+\zeta_3)}$ and extremality of $\lambda_3$ in $\Pi(\lambda_2,\mu_3)$, we obtain that
\begin{align*}
    \eta_{3}=\zeta_{3}=\lambda_{3},\quad\text{and}\quad \eta_{4},\;\zeta_{4}\in\Pi(\lambda_{3},\mu_4).
\end{align*}
Analogously, it can be obtained that
\begin{align*}
    \eta_{N-1}=\zeta_{N-1}=\lambda_{N-1},\quad\text{and}\quad \eta,\;\zeta\in\Pi(\lambda_{N-1},\mu_N),
\end{align*}
and by the extremality of $\lambda$ in $\Pi(\lambda_{N-1},\mu_N)$ and \eqref{QQQ} we end up with  a contradiction. Hence we have the claim.
    \end{proof}
\end{lemma}
Next lemma indicates the relation between  minimizing sets $\mathcal{S}$ and $\mathcal{S}_j$, respectively for \eqref{TAGMMKP} and \eqref{TAGRMKP}.
\begin{lemma}\label{LEM-MINI-SSJ}
 Assume that $\left\{(X_i,\mathcal{B}_{X_i},\mu_i)\right\}_{i=1}^N$ is a collection of probability spaces  such that $X_i$'s are compact subsets of $\mathbb{R}^{n+1}$ and $c:X_1\cdots X_N\to\mathbb{R}$ is a continuous function. For $\mathcal{S}$ and $\mathcal{S}_j$, the minimizing sets for \eqref{TAGMMKP} and \eqref{TAGRMKP}, respectively,  we have
 \begin{align}\label{MINI-SETS}
     \pi_j(\mathcal{S})=\mathcal{S}_j.
 \end{align}
\begin{proof}
        First we note that Proposition \ref{PROP-OPTIMAL} implicitly implies that $\{\phi_i\}_{i=1}^j$ is a solution to the dual of \eqref{TAGRMKP}. 
        Note that by \eqref{PHIONE}  one can consider the $c$-concave potential $\phi_1$, as a $c_j$-concave function as well,
    \begin{align*}
        \phi_1(x)=\inf\left\{c_j(x_1,\ldots,x_j)-\sum_{i=1}^j\phi_i(x_i)\ :\ (x_2,\ldots,x_j)\in\prod_{i=1}^j X_i\right\}.
    \end{align*}
        Additionally, the minimizng sets $\mathcal{S}$ and $\mathcal{S}_j$ are in fact $\partial^c\phi_1$ and $\partial^{c_j}\phi_1$, the $c$- and $c_j$-superdifferential of the potential $\phi_1$,
        \begin{align*}
    &\mathcal{S}=\partial^c\phi_1=\left\{(x_1,\ldots,x_N)\in\prod_{i=1}^NX_i\ :\ c(x_1,\ldots,x_N)=\sum_{i=1}^N\phi_i(x_i)\right\},\\
    &\mathcal{S}_j=\partial^{c_j}\phi_1=\left\{(x_1,\ldots,x_j)\in\prod_{i=1}^j X_i\ :\ c_j(x_1,\ldots,x_j)=\sum_{i=1}^j\phi_i(x_i)\right\}.
\end{align*}
   Let $(x_1,\ldots,x_j)\in \pi_j(\mathcal{S})$. Therefore, there exists $(x_{j+1},\ldots,x_N)\in\prod_{i=j+1}^NX_i$ such that 
   \begin{align}
       (x_1,\ldots,x_j,x_{j+1},\ldots, x_N)\in\mathcal{S}.
   \end{align}
   Moreover, it implies that the infimum in the definition of  $c_j(x_1,\ldots,x_j)$ is attained at $(x_{j+1},\ldots,x_N)$ and we have,
   \begin{align*}
       c_j(x_1,\ldots,x_j)=\sum_{i=1}^j\phi_i(x_i).
   \end{align*}
   Hence $\pi_j(\mathcal{S})\subseteq \mathcal{S}_j$. To show the inverse inclusion let $(x_1,\ldots,x_j)\in \mathcal{S}_j$. Therefore
   \begin{align}\label{WWW1}
       c_j(x_1,\ldots,x_j)=\sum_{i=1}^N\phi_j(x_i).
   \end{align}
   On the other hand, the set $\prod_{i=j+1}^N X_i$ is compact, so by continuity of $c$ and the potentials $\{\phi_i\}_{i=j+1}^N$, for $(x_1,\ldots,x_j)\in \mathcal{S}_j$, the infimum in the definition of $c_j(x_1,\ldots,x_j)$ is attained at some $(x_{j+1},\ldots,x_N)\in\prod_{i=j+1}^NX_i$, and we have
   \begin{align}\label{WWW2}
       c_j(x_1,\ldots,x_j)= c(x_1,\ldots,x_j,x_{j+1},\ldots,x_N)-\sum_{i=j+1}^N\phi_i(x_i).
   \end{align}
By \eqref{WWW1} and \eqref{WWW2} we have $(x_1,\ldots,x_j,x_{j+1},\ldots,x_N)\in\mathcal{S}$. Therefore, $(x_1,\ldots,x_j)\in\pi_j(\mathcal{S})$.

\end{proof}
\end{lemma}
As an immediate consequence of Theorem \ref{TH-CP-EXT}, Proposition \ref{PROP-OPTIMAL}, and Lemmas \ref{LEM-TRI-EXT} and \ref{LEM-MINI-SSJ}, we have the following applicable theorem.
\begin{theorem}\label{TH-UNI-MULT-CP-EXT}
Let $\left\{(X_i,\mathcal{B}_{X_i},\mu_i)\right\}_{i=1}^N$ be a family of probability spaces  where $X_i$'s are compact subsets of $\mathbb{R}^n$ and $P_i=\{X_{ik}\}_{k=1}^{K_i}$ be a Borel ordered partition of $X_i$, for $i=2,\ldots,N$. Assume that $c:X_1\times\cdots\times X_N\to\mathbb{R}$ and $\mathcal{S}$ is a minimizing set for \eqref{TAGMMKP}. Then if $\mathcal{S}$ is $(c,P_N)$-extreme and $\pi_j(\mathcal{S})$ is $(c_j,P_j)$-extreme for $j=2,\ldots,N-1$, then the solution to \eqref{TAGMMKP} is unique.
\begin{proof}
By Lemma \ref{LEM-MINI-SSJ} we know that $\pi_j(\mathcal{S})$ is a minimizing set for \eqref{TAGRMKP}, for $j=2,\ldots,N-1$. By assumption, the set $\mathcal{S}$ is $(c,P_N)$-extreme. Thus, via Theorem \ref{TH-CP-EXT}, this implies that each optimal plan $\lambda$ of \eqref{TAGMMKP} is an extreme point of $\Pi(\lambda_{N-1},\mu_N)$. On top of that, by Proposition \ref{PROP-OPTIMAL}, $\lambda_j=\pi_j\#\lambda$ is an optimal plans for \eqref{TAGRMKP}, for $j=2,\ldots, N-1$. But by assumptions, $\pi_j(\mathcal{S})$'s are $(c_j,P_j)$-extreme. Therefore, again by Theorem \ref{TH-CP-EXT}, the optimal plan $\lambda_j$ is an extreme point of $\Pi(\lambda_{j-1},\mu_j)$. By applying Lemma \ref{LEM-TRI-EXT}, we obtain that $\lambda$ is an extreme point of $\Pi(\mu_1,\ldots,\mu_N)$. Now, since each optimal plan of \eqref{TAGMMKP} is an extreme point of $\Pi(\mu_1,\ldots,\mu_N)$, then the uniqueness follows.
    
\end{proof}
\end{theorem}

\begin{remark}\label{REM-TH-UI-MULT-CP-EXT}
    The same result holds if we replace $(c,P_N)$- and $(c_j,P_j)$-extremality assumptions on $\mathcal{S}$  and $\pi_j(\mathcal{S})$ by $c$- and $c_j$-extremality, when one works with the partitions $P_i=\{X_i\}$.
\end{remark}

\end{document}